\documentclass{article}
\usepackage{amssymb,amsmath,amsthm}
\usepackage{graphicx}

\newtheorem{theorem}{Theorem}
\newtheorem{lemma}[theorem]{Lemma}

\newtheorem{corollary}[theorem]{Corollary}

\theoremstyle{definition}
\newtheorem{example}{Example}
\newtheorem{definition}{Definition}
\newtheorem{remark}{Remark}

\date{}

\title{\Large \textbf{Niebrzydowski Algebras and Trivalent Spatial Graphs}}

\author{
Paige Graves\footnote{Paige.Graves@laverne.edu}\and
Sam Nelson \footnote{Sam.Nelson@cmc.edu. Partially supported 
by Simons Foundation Collaboration Grant 316709.} \and 
Sherilyn Tamagawa\footnote{Tamagawa@math.ucsb.edu}}

\begin{document}
\maketitle

\begin{abstract}
We introduce \textit{Niebrzydowski algebras}, algebraic structures with a 
ternary operation and a partially defined multiplication, with axioms 
motivated by the Reidemeister moves for $Y$-oriented trivalent spatial 
graphs and handlebody-links. As part of this definition, we identify 
generating sets of $Y$-oriented Reidemeister moves. We give some examples to 
demonstrate that the counting invariant can distinguish some $Y$-oriented
trivalent spatial graphs and handlebody-links.
\end{abstract}

\parbox{5.5in}{\textsc{Keywords:} Biquasiles, tribrackets, spatial graphs, 
handlebody-links

\smallskip

}

\section{\large \textbf{Introduction}}

In \cite{MN} an algebraic structure known as \textit{ternary 
quasigroup} was used for coloring planar knot complements.
In \cite{dsn2}, a related algebraic structure known as \textit{biquasile} was 
introduced with axioms derived from the Reidemeister moves for oriented
knots and links with elements used to color vertices for \textit{dual 
graph diagrams} corresponding to regions in the planar complement of a
knot diagram. Biquasile coloring invariants were enhanced with Boltzmann 
weights in \cite{cnn} which are related to cocycles in the cohomology
theory for ternary quasigroups in \cite{MN2}. Biquasile colorings and
ternary quasigroups were used to distinguish orientable surface links
in \cite{KN} and \cite{MN}. In \cite{NP}
structures known as \textit{virtual tribrackets} were defined for
coloring planar complements of oriented virtual link diagrams. 

\textit{Spatial graphs}, i.e. knotted graphs in $\mathbb{R}^3$, were considered
in terms of diagrammatic moves in \cite{K}. Handlebody-links and their
Reidemeister moves were considered in \cite{I}, and \textit{Y-orientations}
of trivalent spatial graphs and handlebody-links were considered in \cite{II}.
Quandle and biquandle colorings of handlebody-links have been considered
in \cite{IN, L}, including algebraic structures such as \textit{G-families of 
quandles and biquandles}, \textit{partially multiplicative quandles and 
biquandles} and \textit{qualgebras}
In \cite{NZ}, the algebraic structure known as \textit{bikei} (see \cite{AN})
was extended by adding a partially defined multiplication to define a 
coloring structure for unoriented twisted virtual handlebody-links called 
\textit{twisted virtual bikeigebras}.

In this paper we extend tribrackets to define
coloring structures for planar complements of $Y$-oriented trivalent spatial 
graphs and $Y$-oriented handlebody-links, obtaining objects we call 
\textit{Niebrzydowski Algebras}. These are sets with a ternary operation $[,,]$
called a \textit{Niebrzydowski tribracket} and a partially defined
multiplication operation
satisfying axioms coming from the Reidemeister moves for oriented
trivalent spatial graphs. A certain subset of these objects satisfy the
additional condition required by the IH-move, making them suitable for
coloring $Y$-oriented handlebody-links. 

The paper is organized as follows. In Section \ref{T} we review the 
Reidemeister moves for $Y$-oriented trivalent spatial graphs and 
handlebody-links, identifying our preferred generating set of
oriented moves. In Section \ref{B} we introduce our algebraic structures
and coloring rules and show that the cardinality of the set of colorings 
of a diagram 
representing an oriented trivalent spatial graph or a $Y$-oriented 
handlebody-link defines an integer-valued invariant. We compute some 
explicit examples of the new invariant. We conclude in Section \ref{Q}
with some questions for future work.

\section{\large \textbf{Trivalent Spatial Graphs}}\label{T}

A (non-oriented) trivalent spatial graph is a finite graph with all vertices degree 3 embedded in $\mathbb R^3$, where two embedded graphs are considered equivalent if there is an ambient isotopy $\mathbb R^3 \rightarrow \mathbb R^3$ taking one to the other.

From \cite{K}
we know that two (non-oriented) trivalent spatial graphs are equivalent iff their diagrams are related by a finite sequence of the Reidemeister moves

\begin{center}
\includegraphics[height=0.8in]{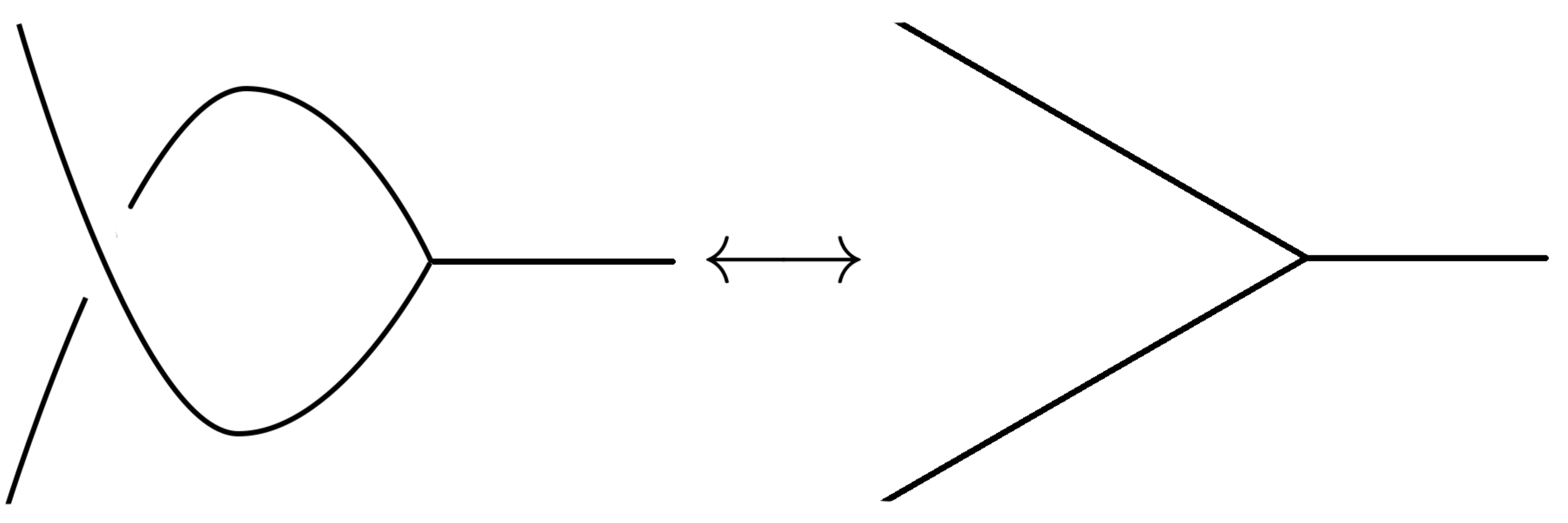} \phantom{HELLO!!}
\includegraphics[height=0.8in]{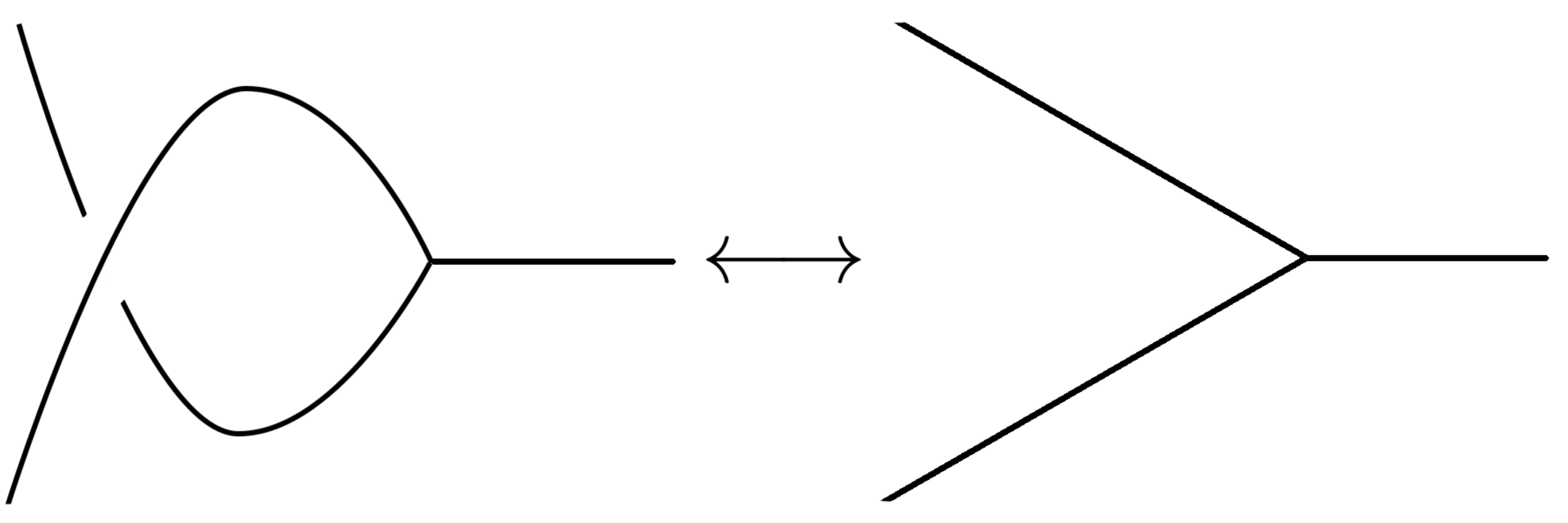}

\includegraphics[height=0.8in]{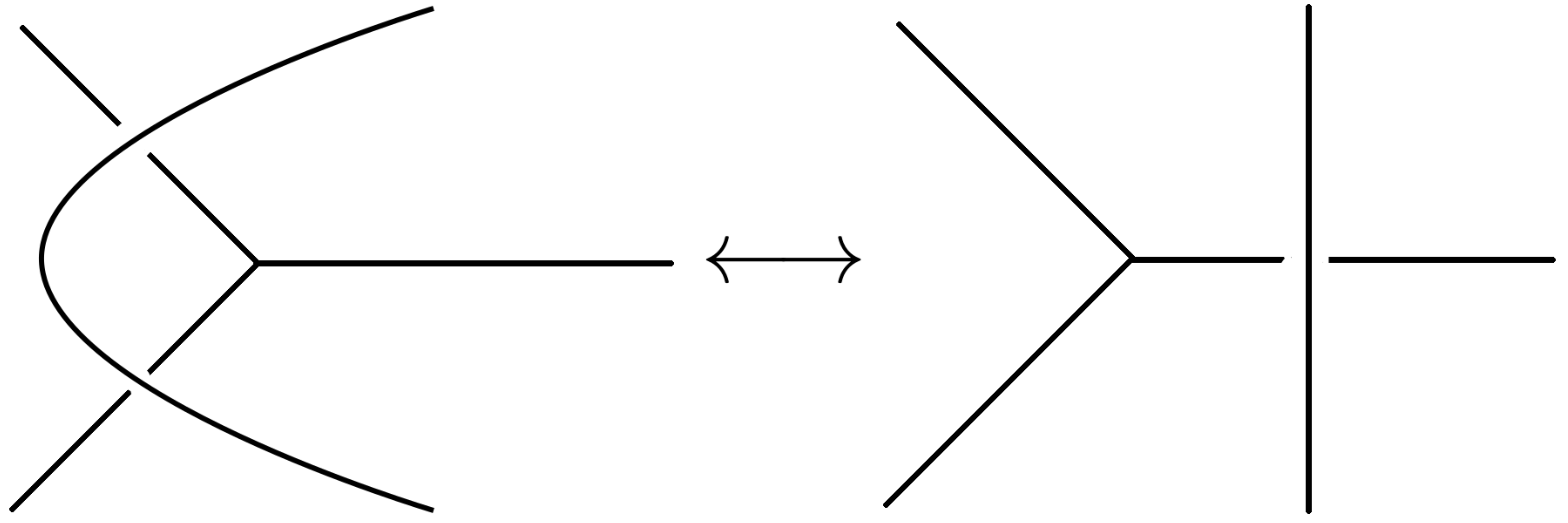} \phantom{HELLO!!}
\includegraphics[height=0.8in]{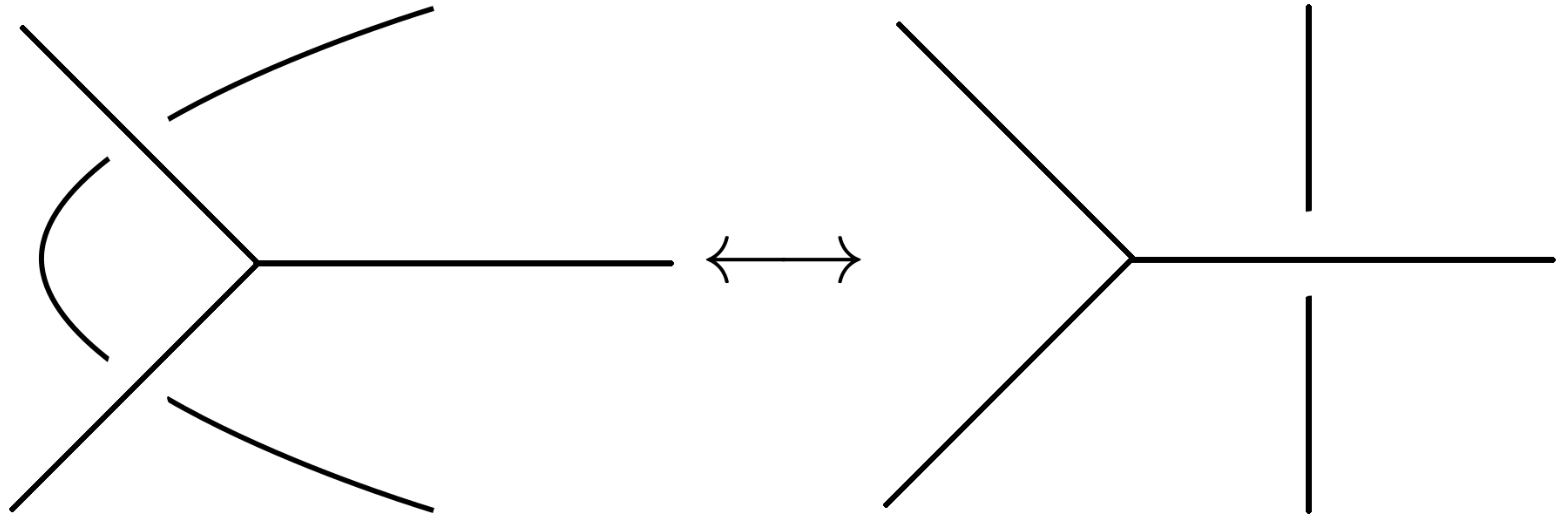}
\end{center}
and the classical Reidemeister moves.

We orient the edges of a trivalent spatial graph, forbidding sources and sinks at the vertices. Such an orientation is known as a \textit{$Y$-orientation}.

Then, by listing all the possible orientations on all edges of each of $R4$ and $R5$, we see that we have the  generating list of moves below.

\begin{center}
\begin{tabular}{cccc}
$R4.1$ &$ R4.2$ & $R4.3$ & $R4.4$ \\
\includegraphics[height=0.4in]{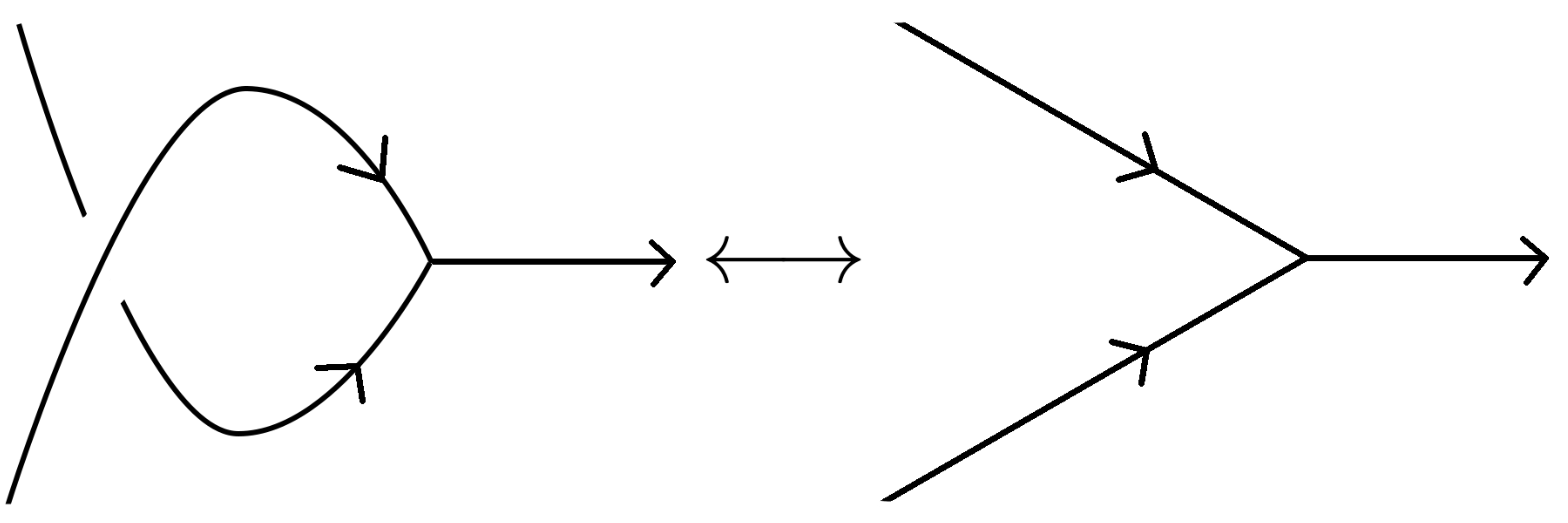} &
\includegraphics[height=0.4in]{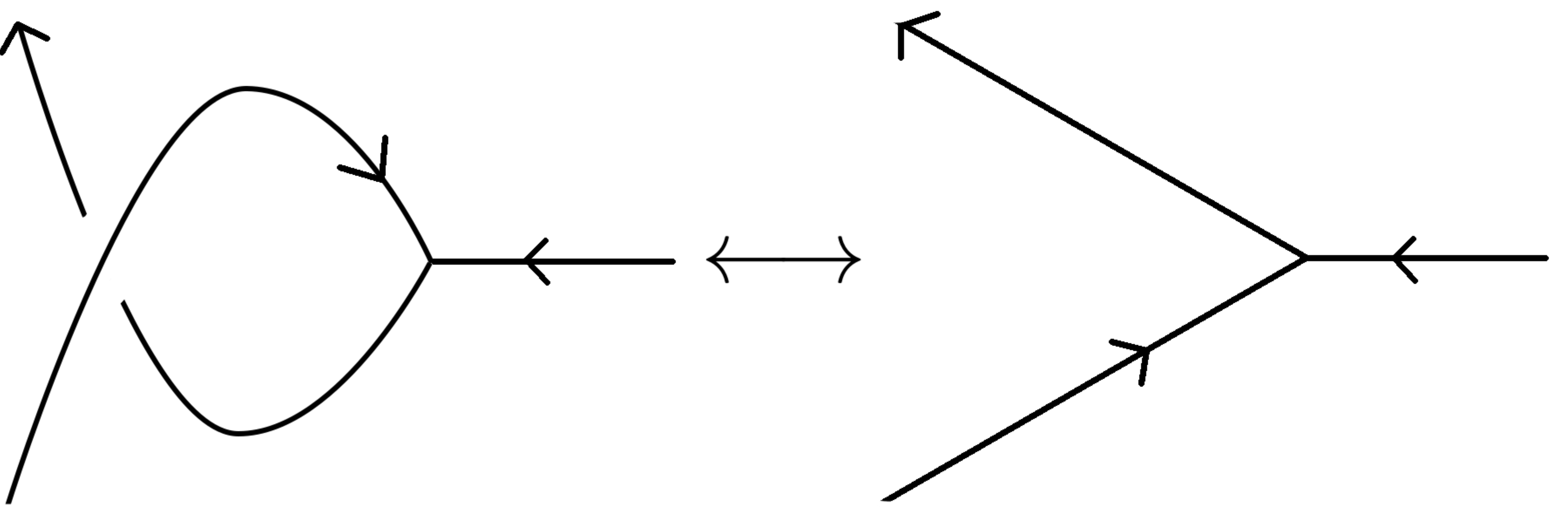}&
\includegraphics[height=0.4in]{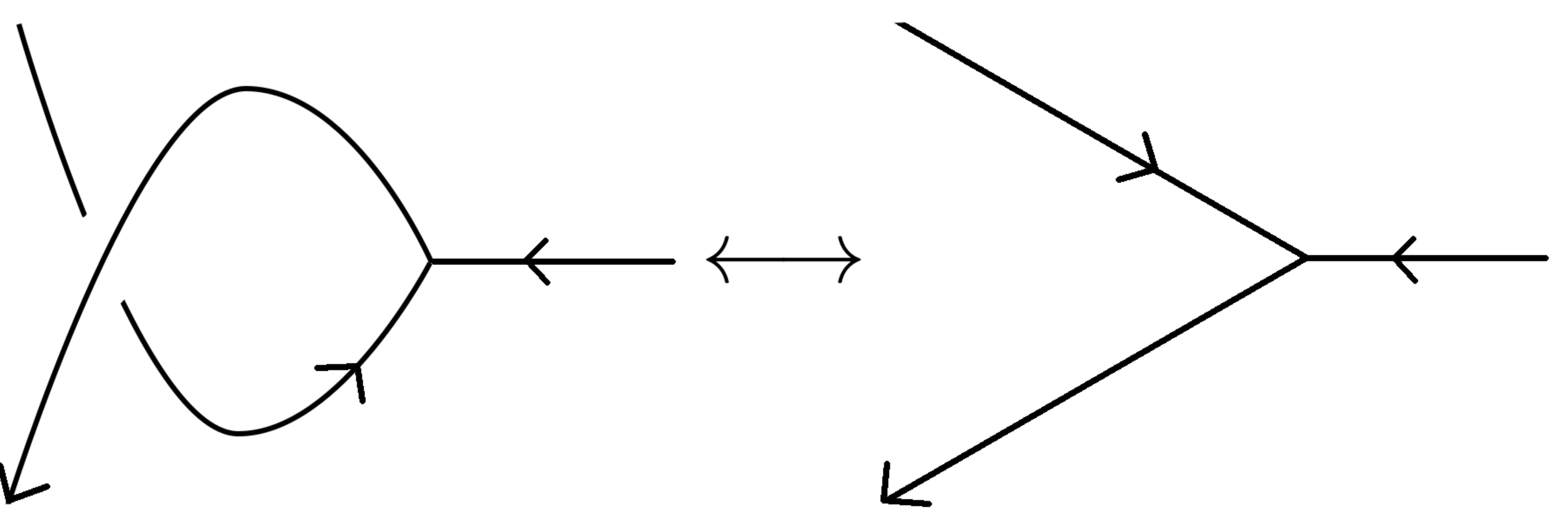}&
\includegraphics[height=0.4in]{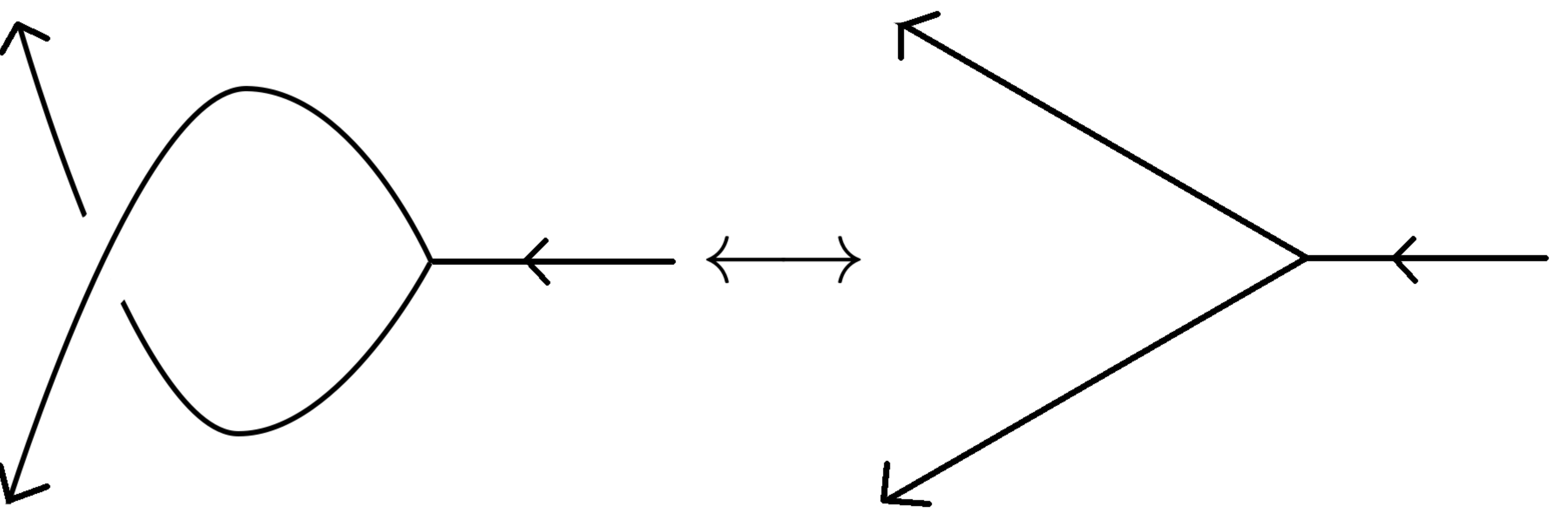}\\ 
\vspace{1em} \\
$R4.5$ & $R4.6$ & $R4.7$ & $R4.8$ \\
\includegraphics[height=0.4in]{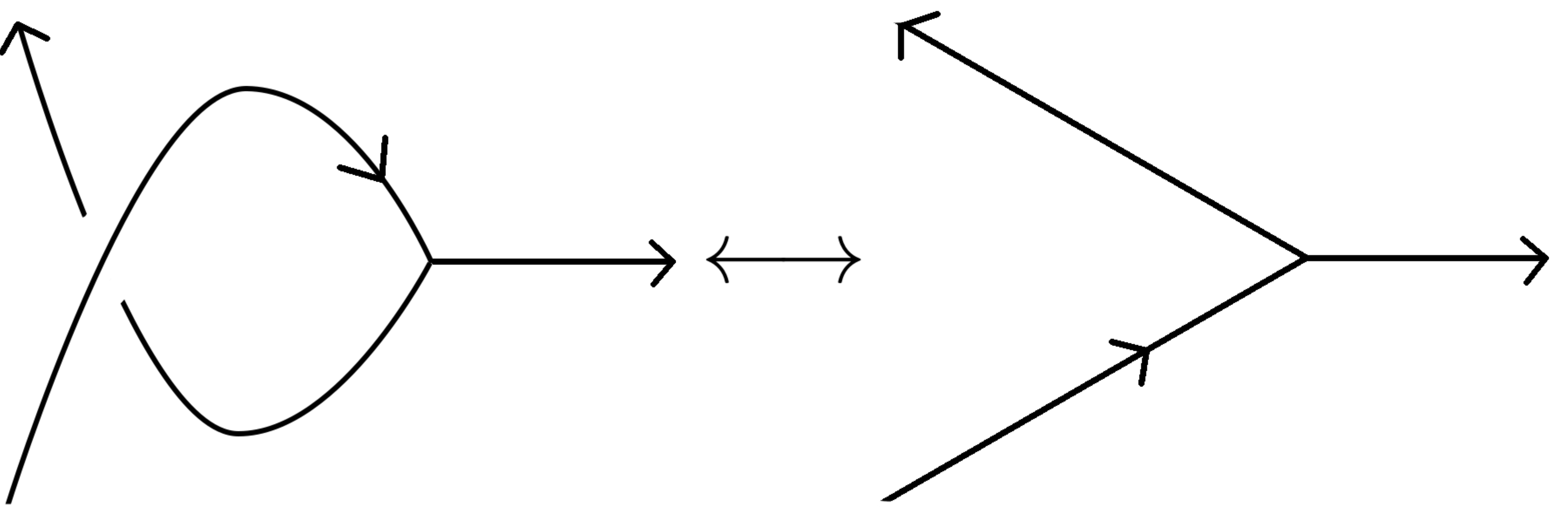} &
\includegraphics[height=0.4in]{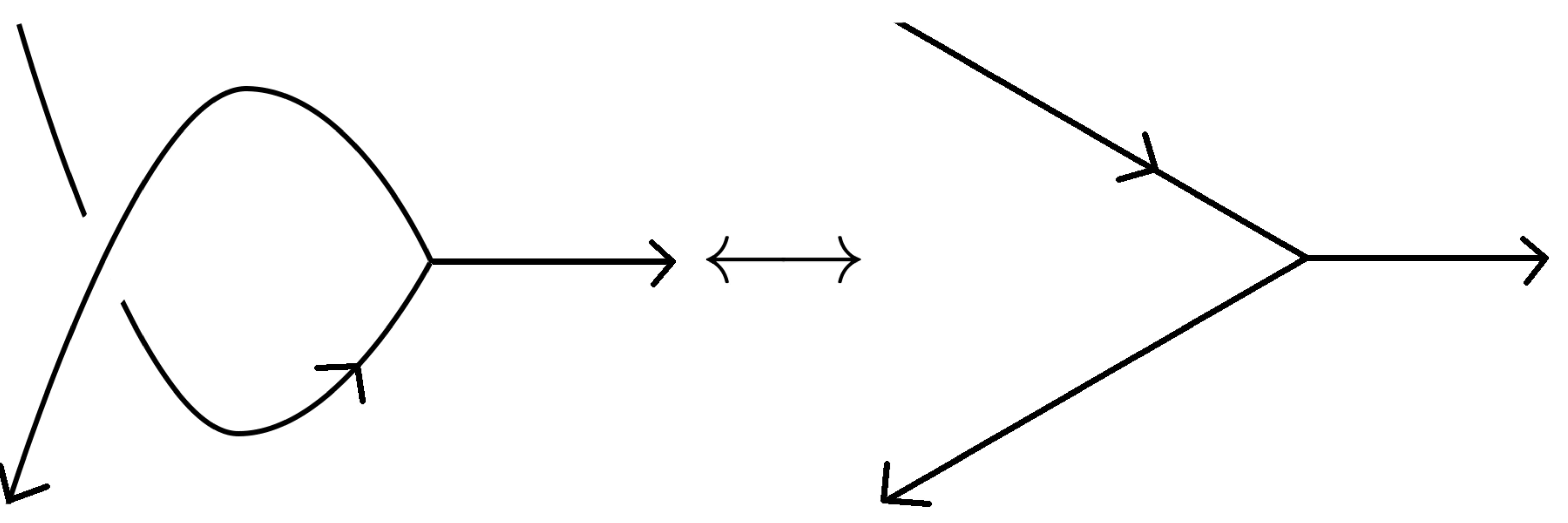} & 
\includegraphics[height=0.4in]{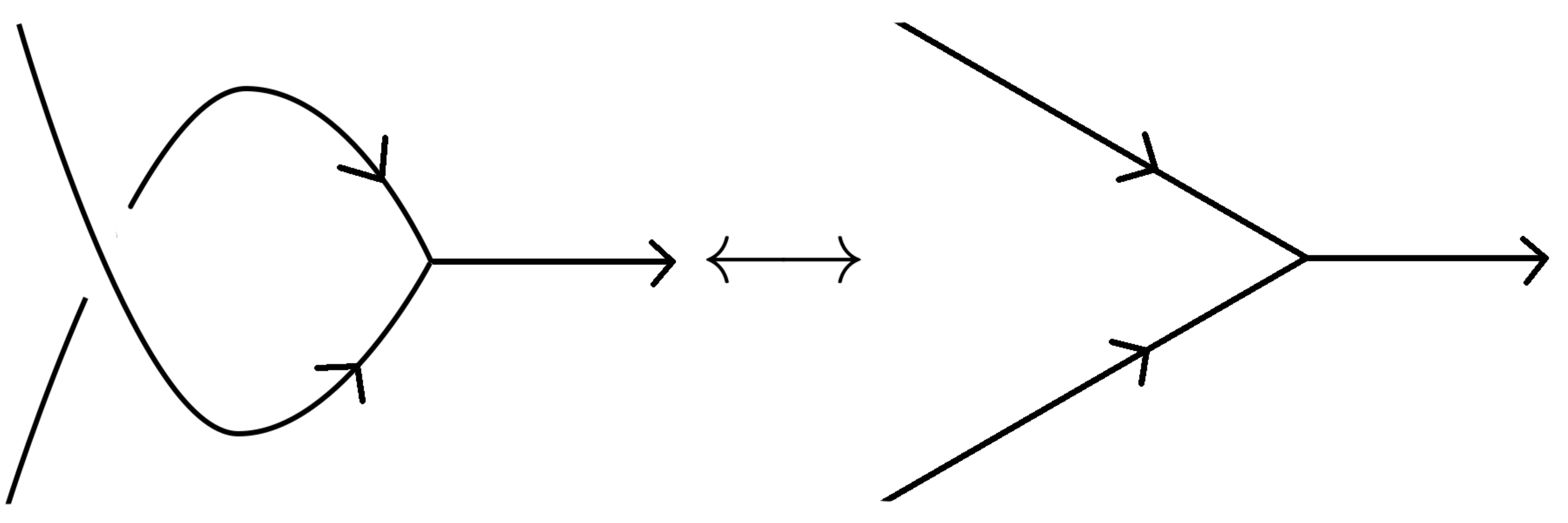}&
\includegraphics[height=0.4in]{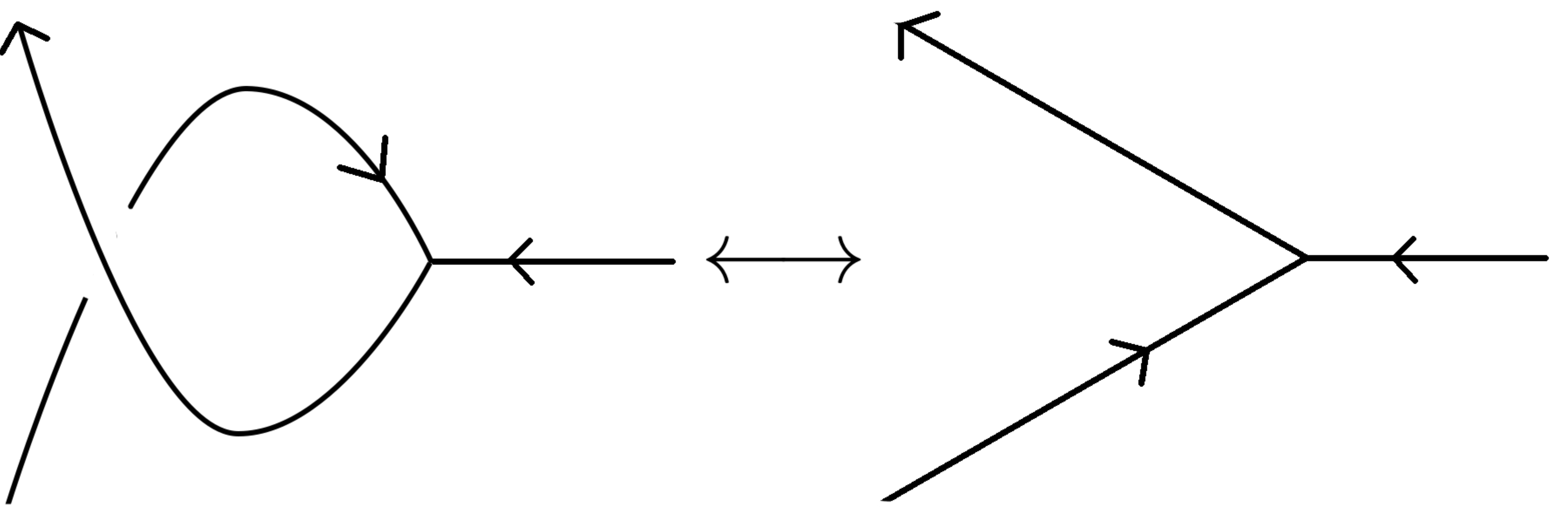} \\
\vspace{1em} \\
$R4.9$ & $R4.10$ & $R4.11$ & $R4.12$ \\
\includegraphics[height=0.4in]{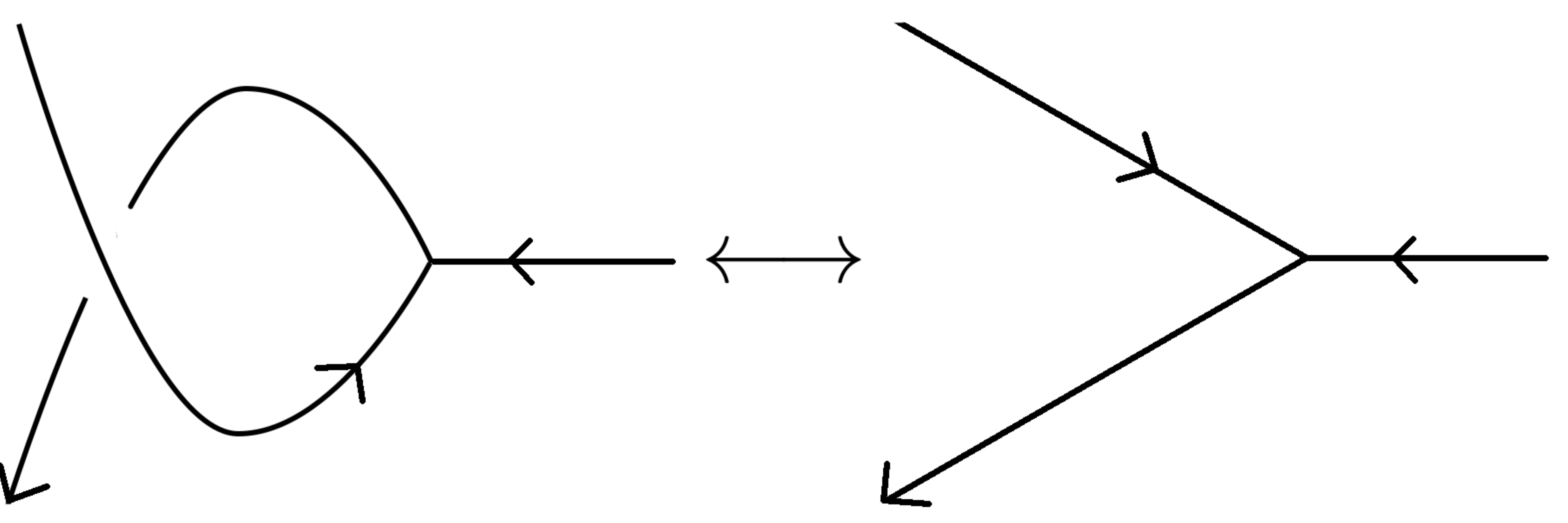} & 
\includegraphics[height=0.4in]{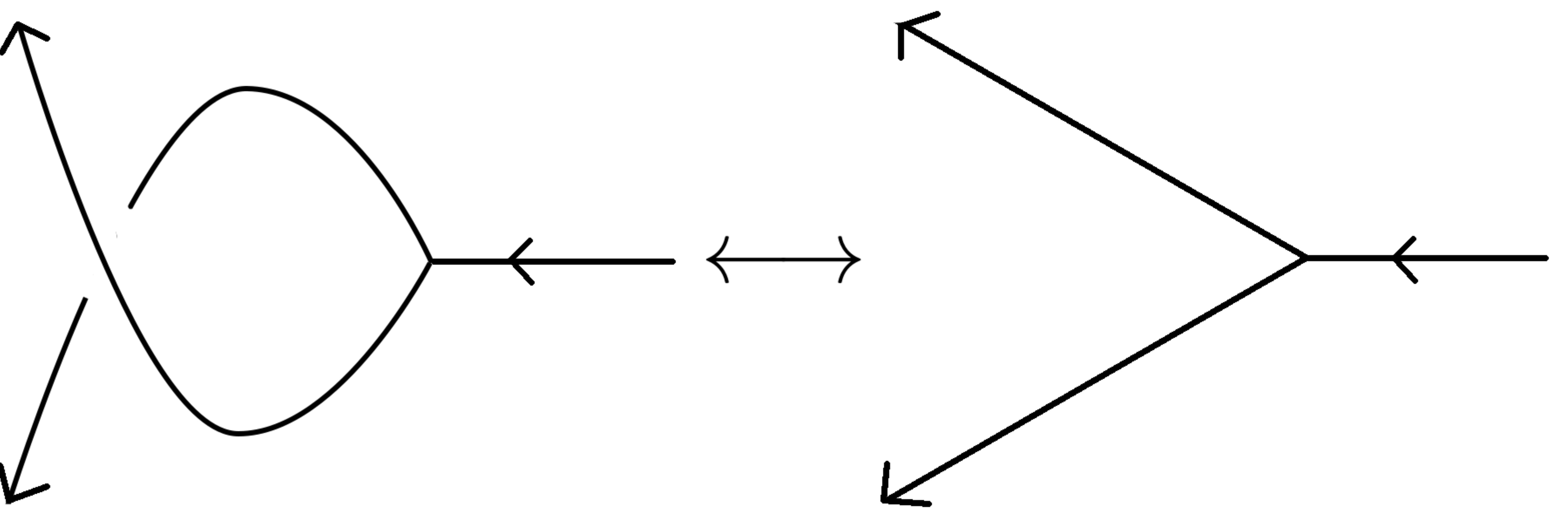}&
\includegraphics[height=0.4in]{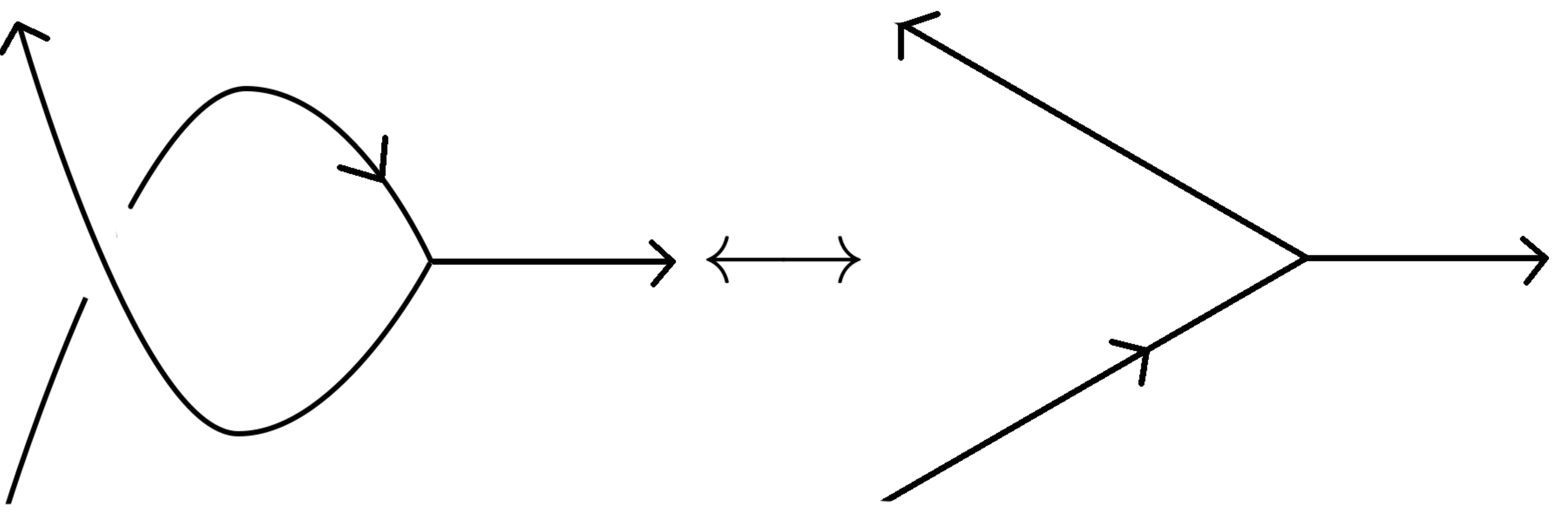} & 
\includegraphics[height=0.4in]{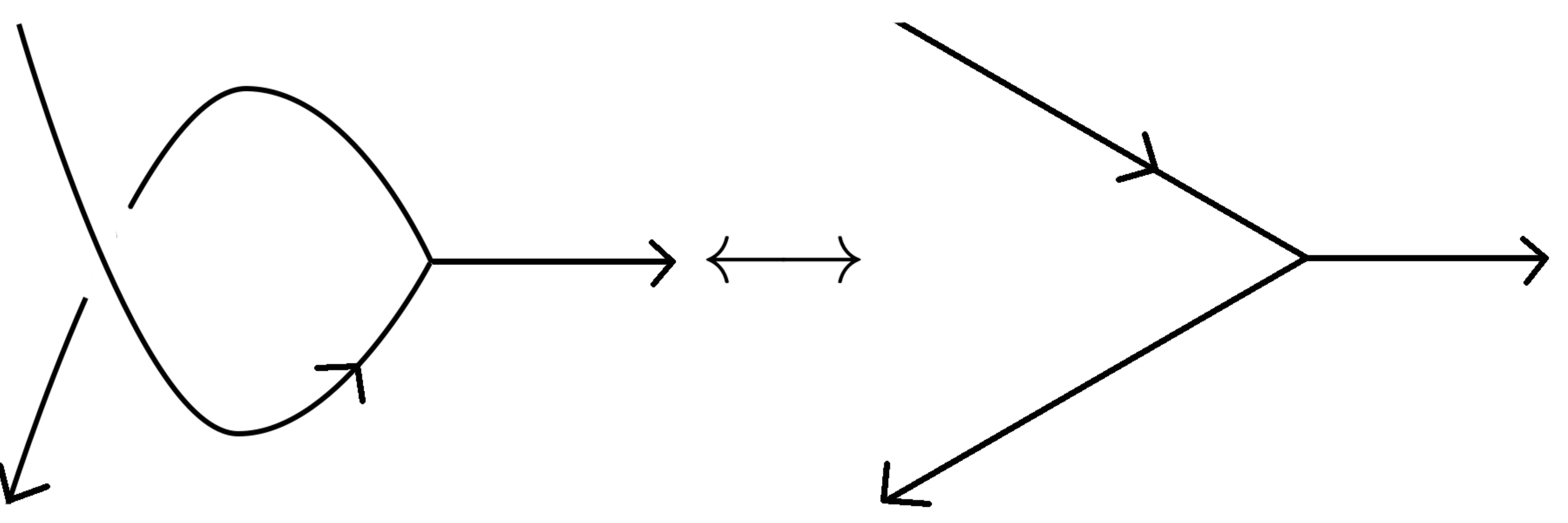} \\
\vspace{1em} \\
$R5.1$ & $R5.2$ & $R5.3$ & $R5.4$ \\
\includegraphics[height=0.4in]{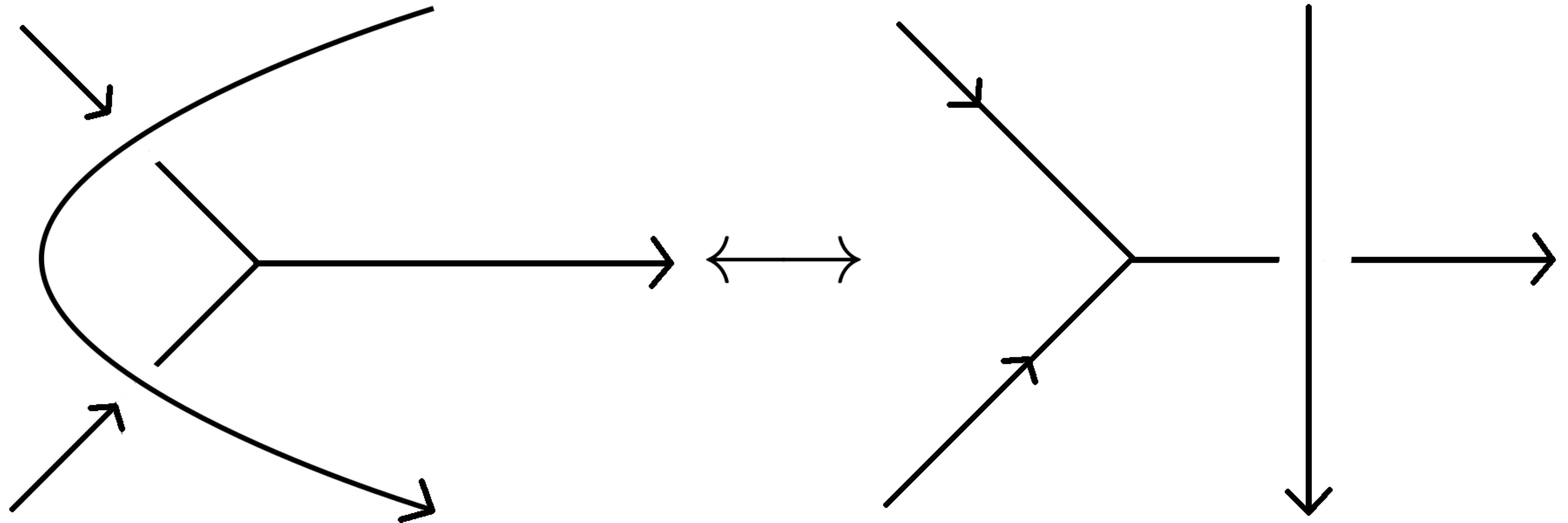} & 
\includegraphics[height=0.4in]{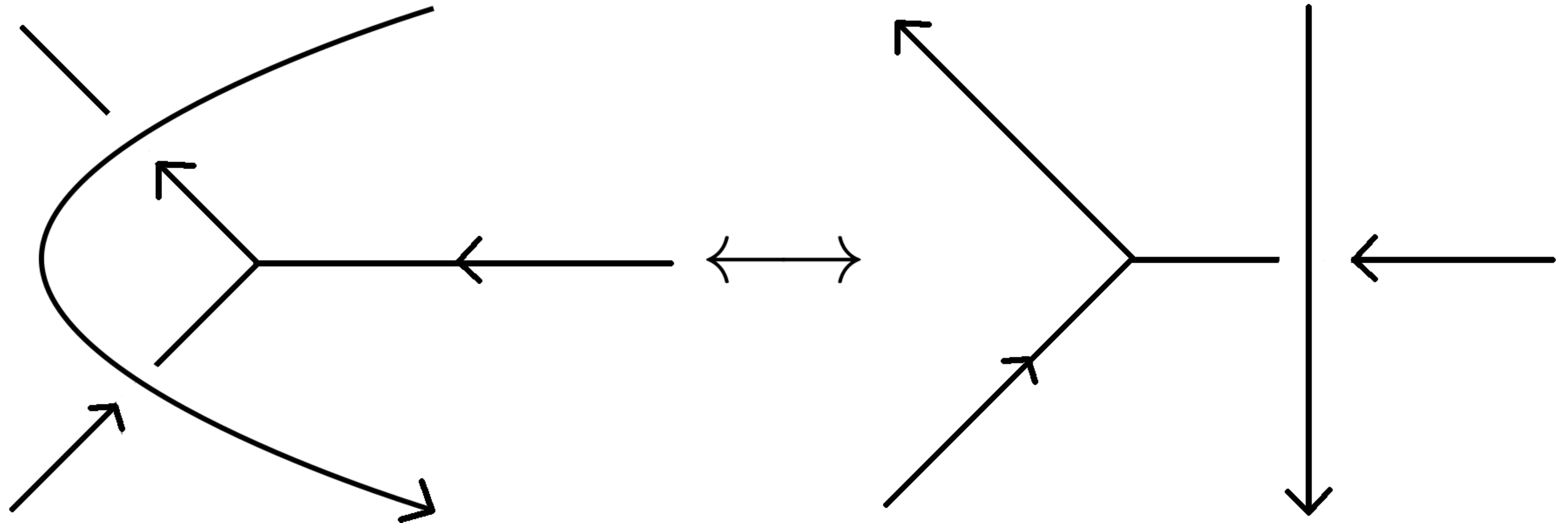}&
\includegraphics[height=0.4in]{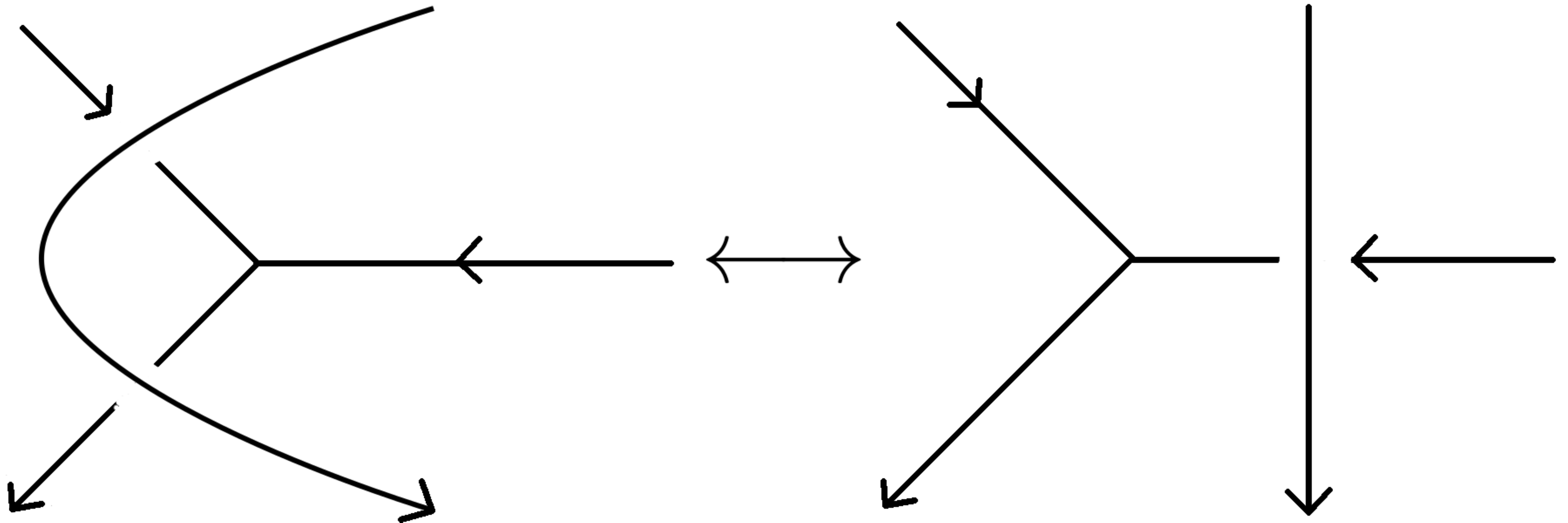} & 
\includegraphics[height=0.4in]{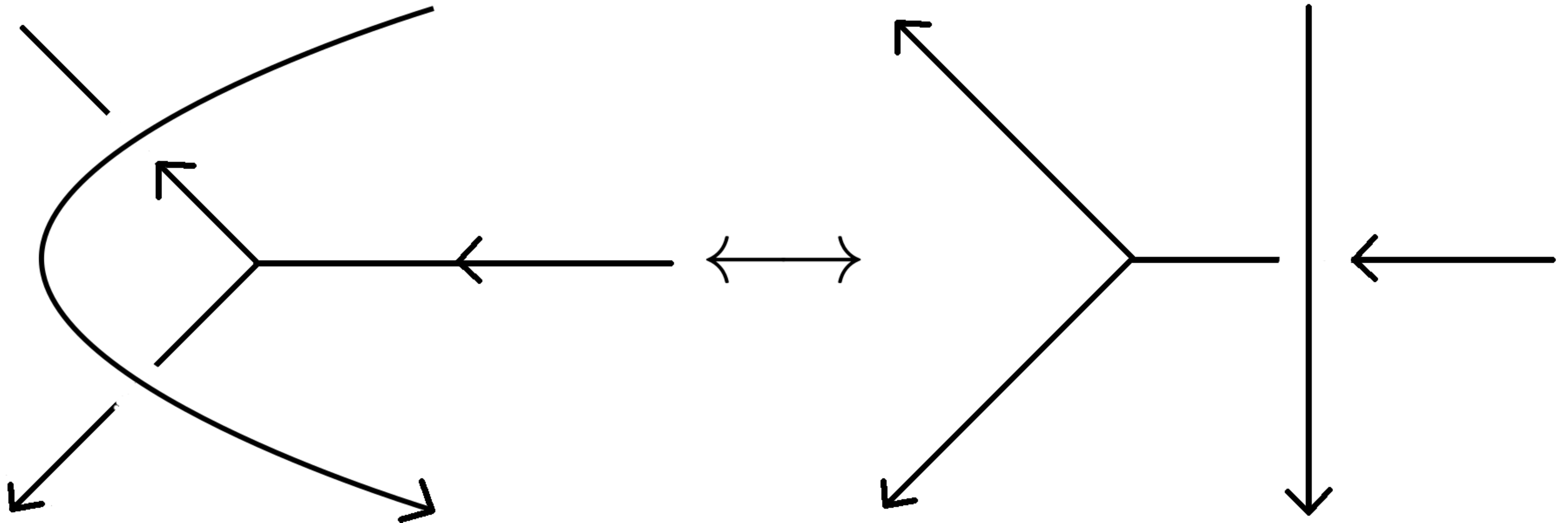} \\
\end{tabular}\end{center}\begin{center}\begin{tabular}{cccc}
$R5.5$ & $R5.6$ & $R5.7$ & $R5.8$ \\
\includegraphics[height=0.4in]{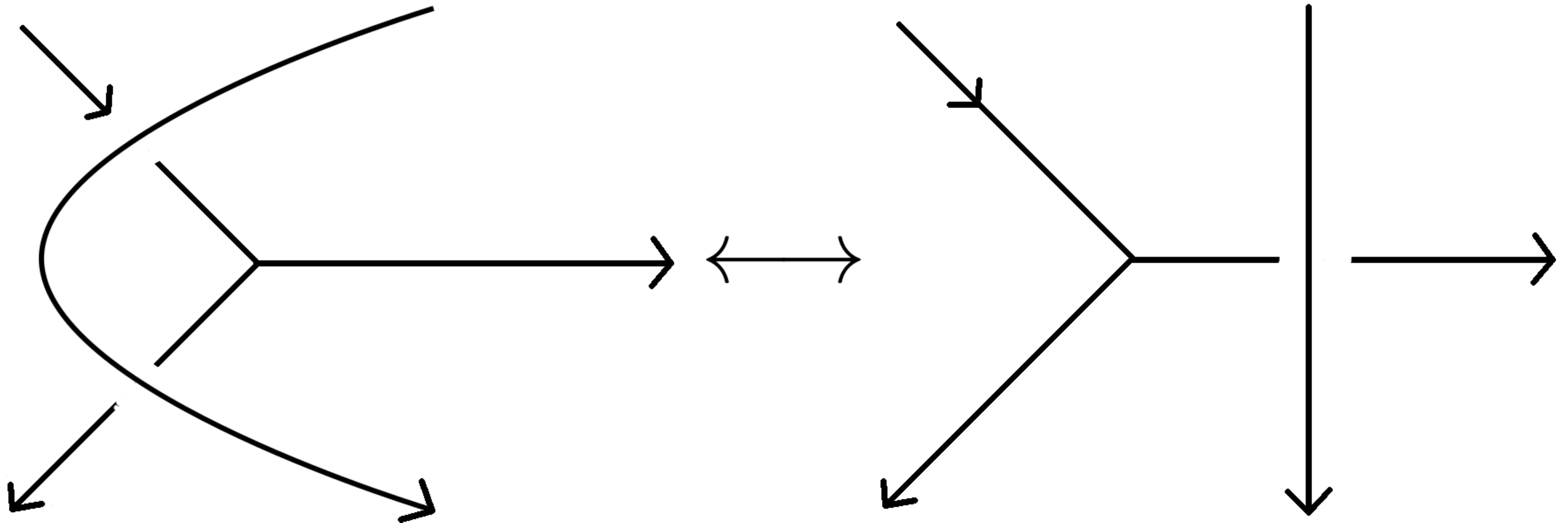} &
\includegraphics[height=0.4in]{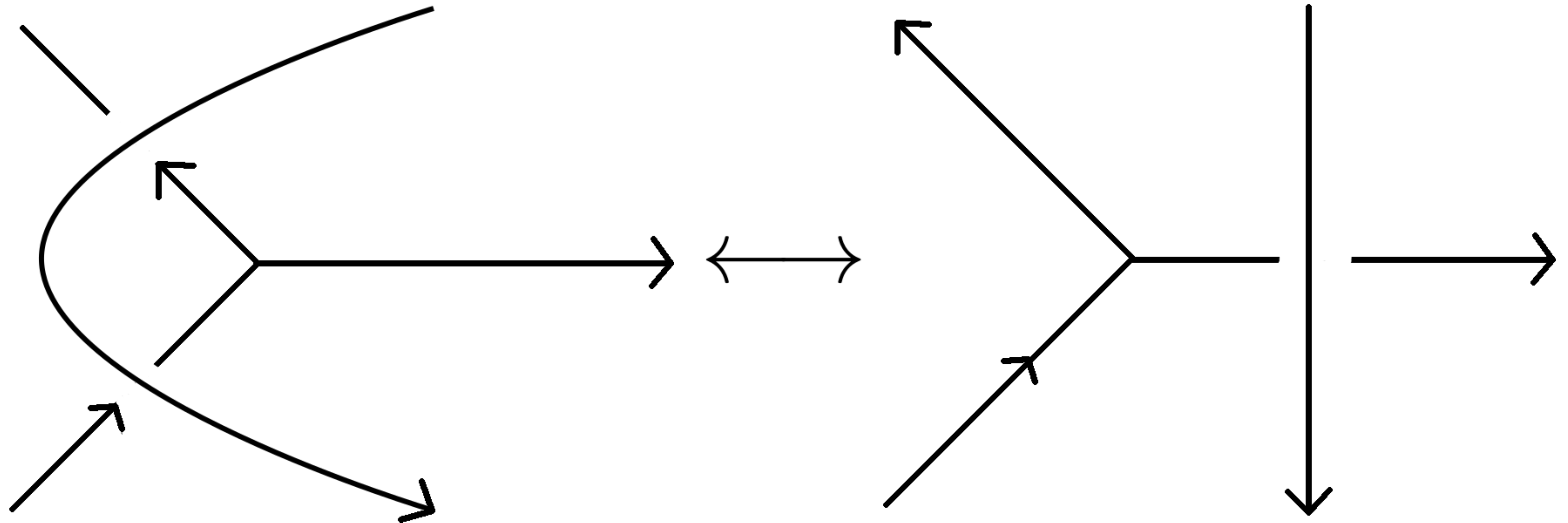} & 
\includegraphics[height=0.4in]{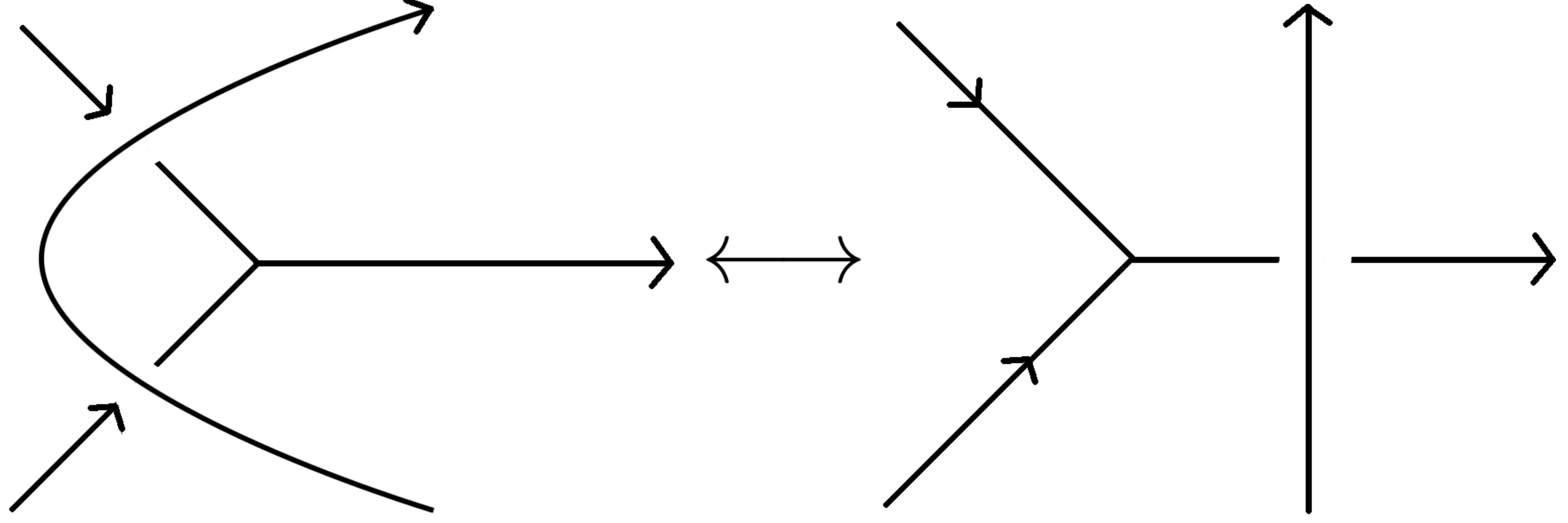} & 
\includegraphics[height=0.4in]{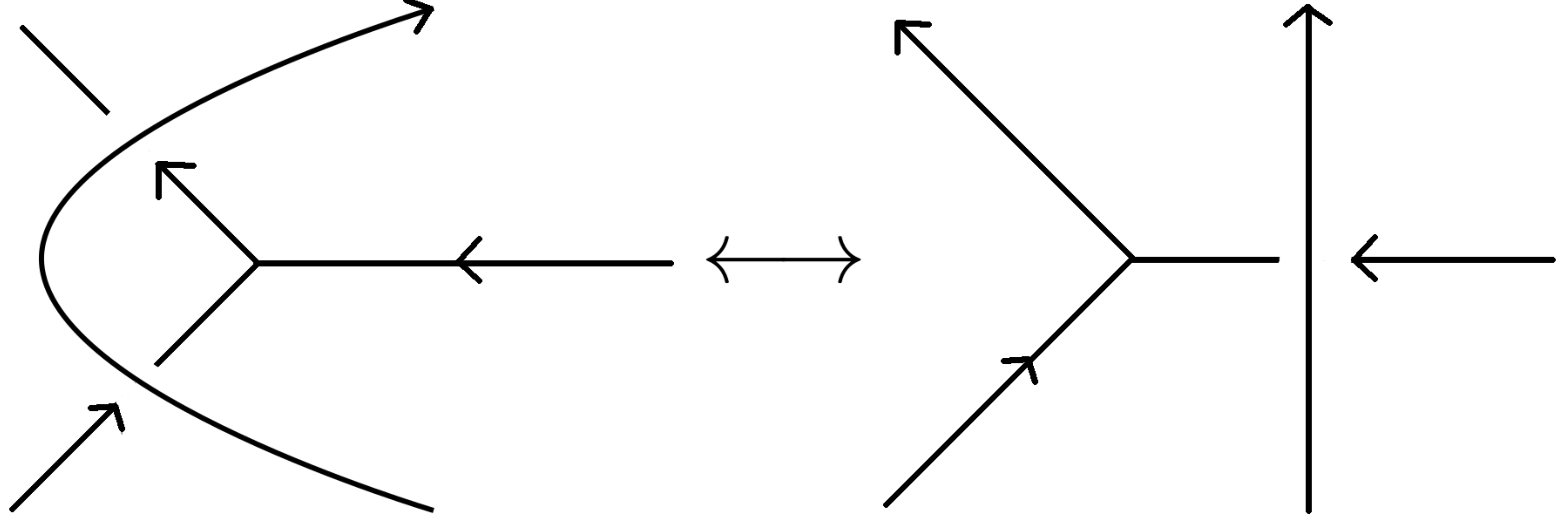} \\
\vspace{1em} \\
$R5.9$ & $R5.10$ & $R5.11$ & $R5.12$ \\
\includegraphics[height=0.4in]{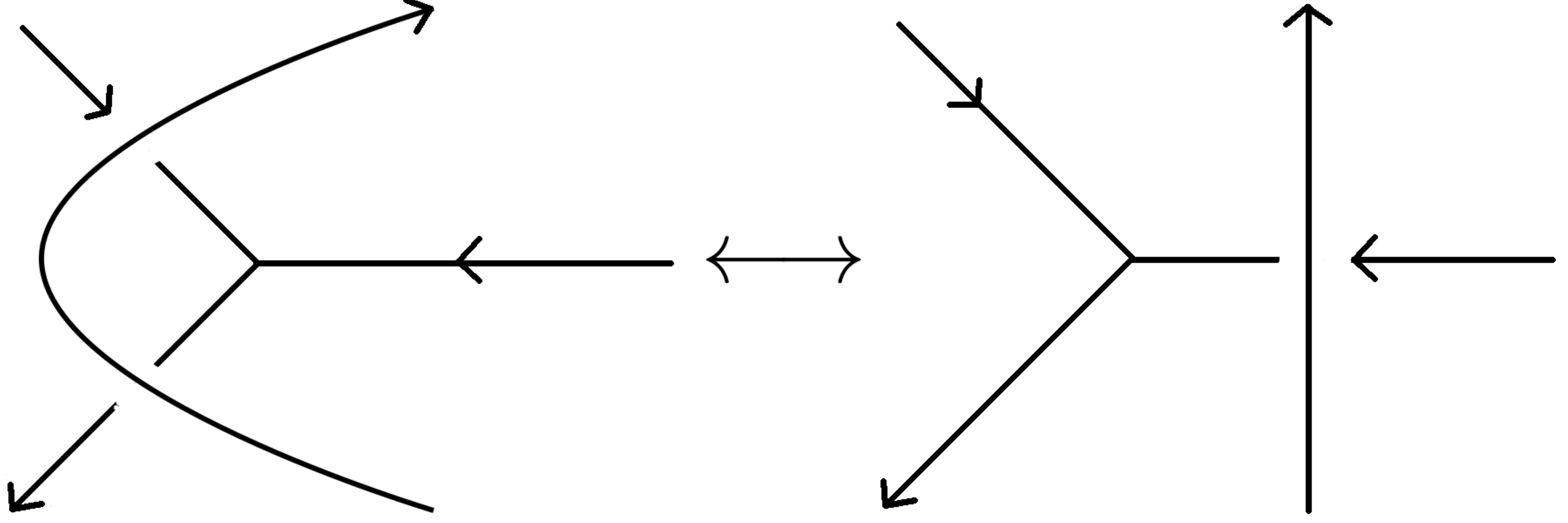} & 
\includegraphics[height=0.4in]{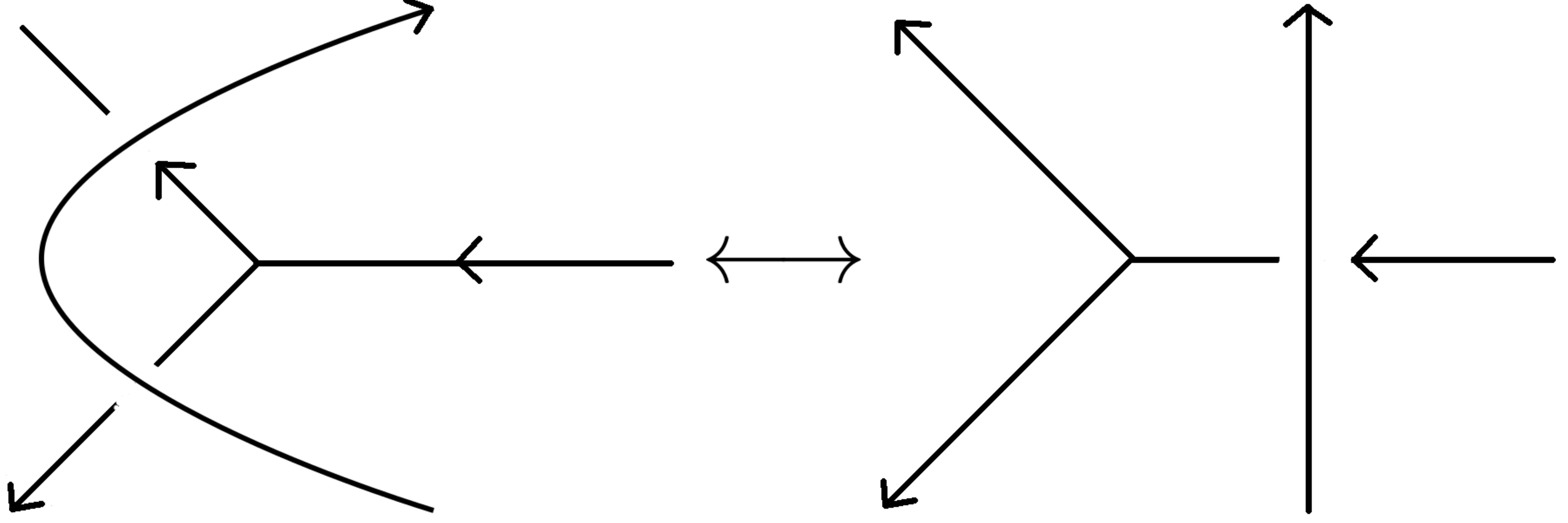} &
\includegraphics[height=0.4in]{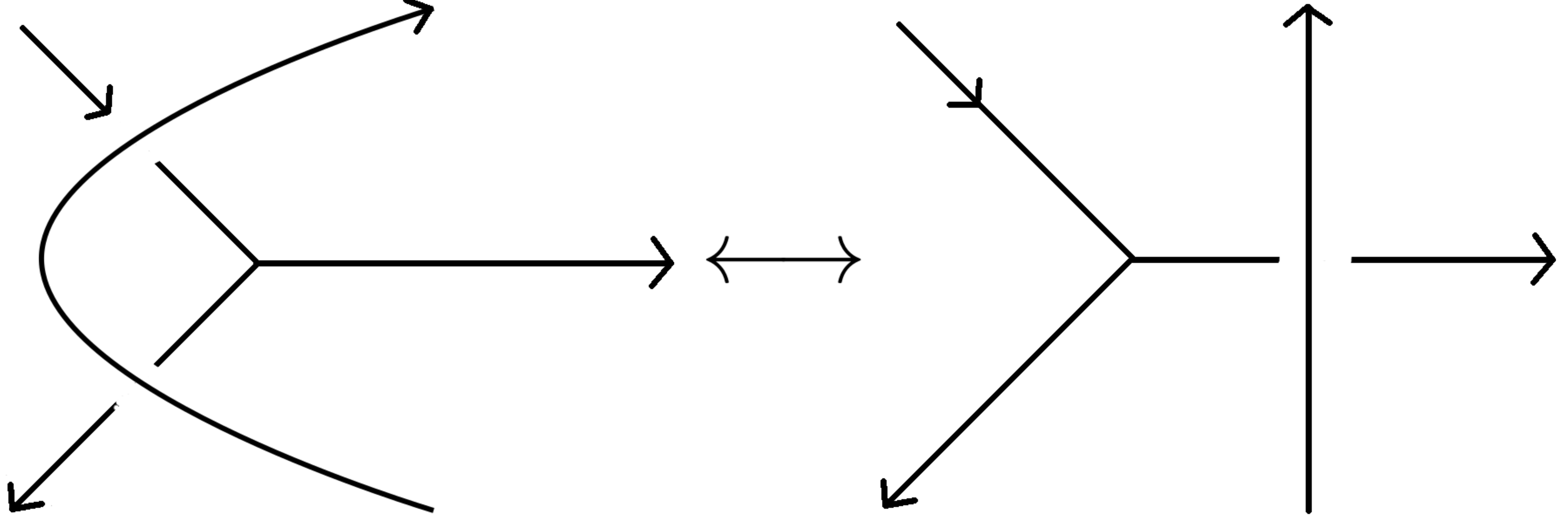} & 
\includegraphics[height=0.4in]{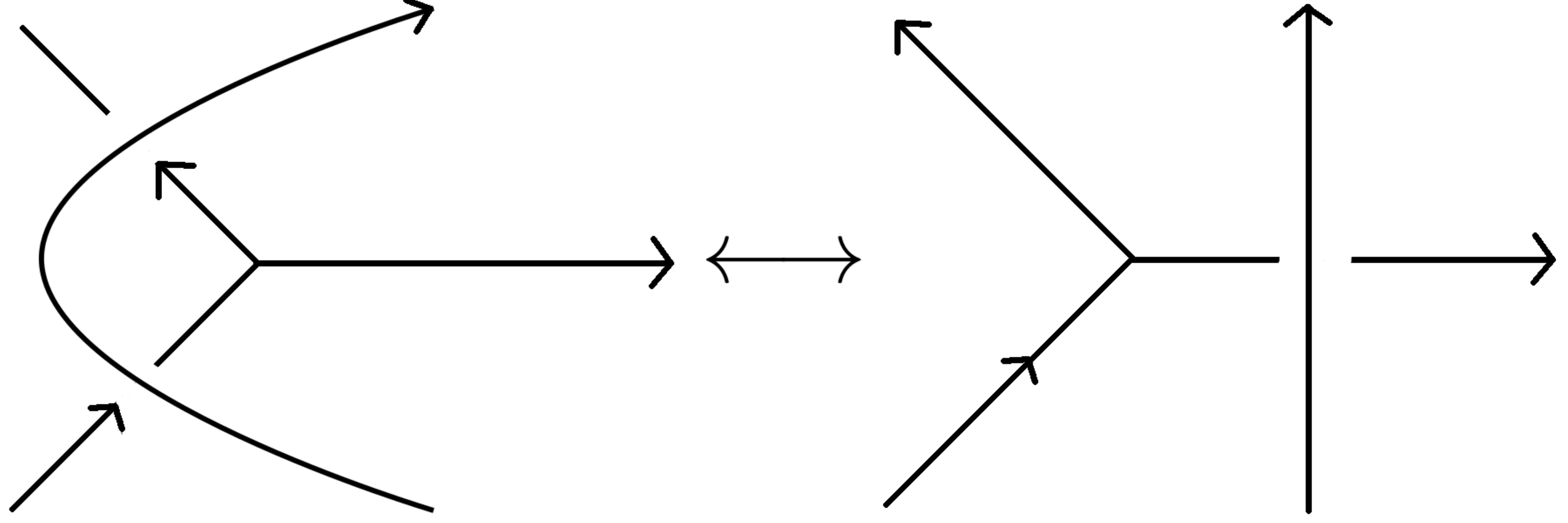} \\
\vspace{1em}\\
$R5.13$ & $R5.14$ & $R5.15$ & $R5.16$ \\
\includegraphics[height=0.4in]{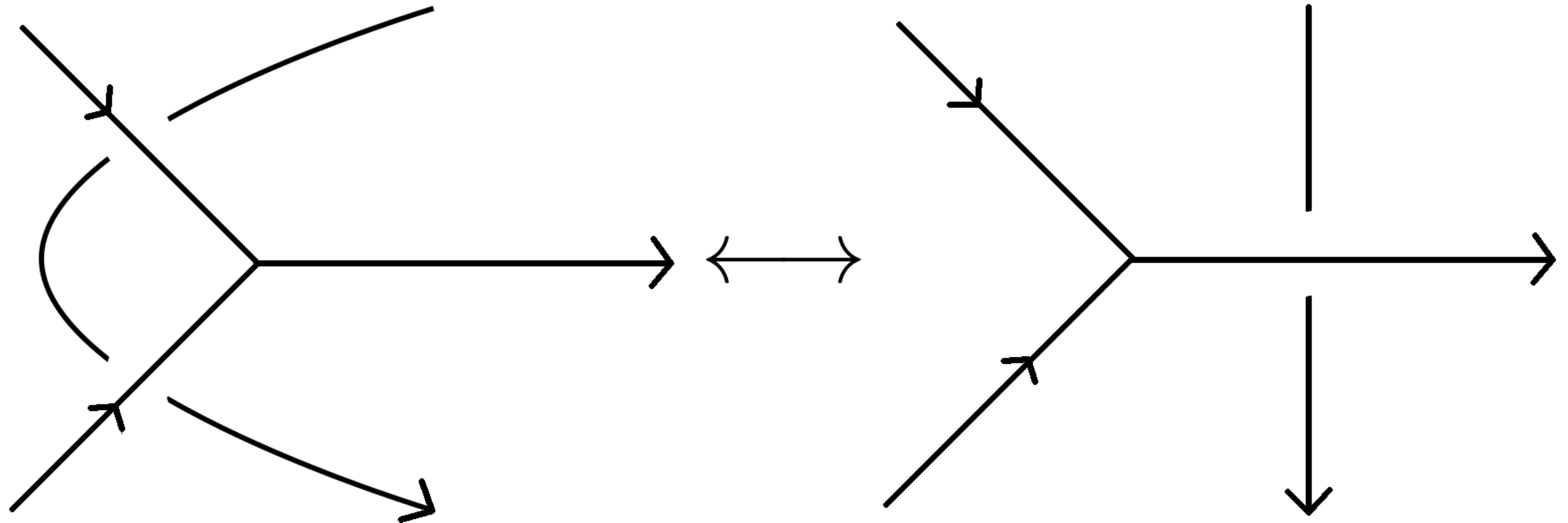} &
\includegraphics[height=0.4in]{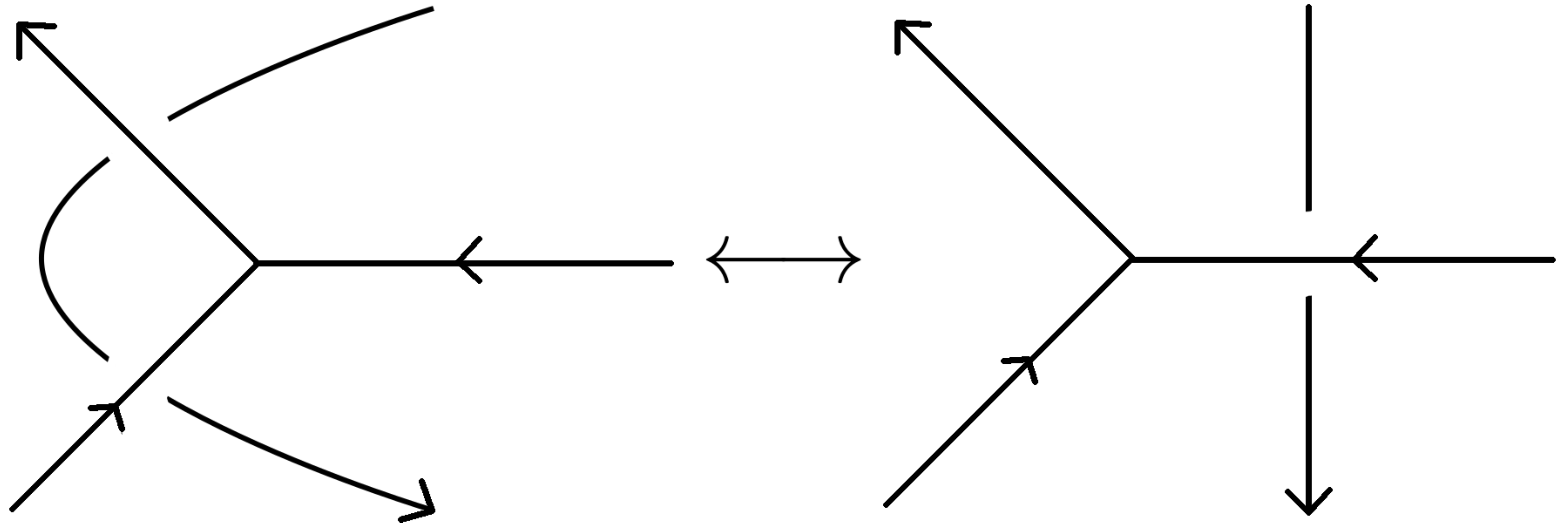} & 
\includegraphics[height=0.4in]{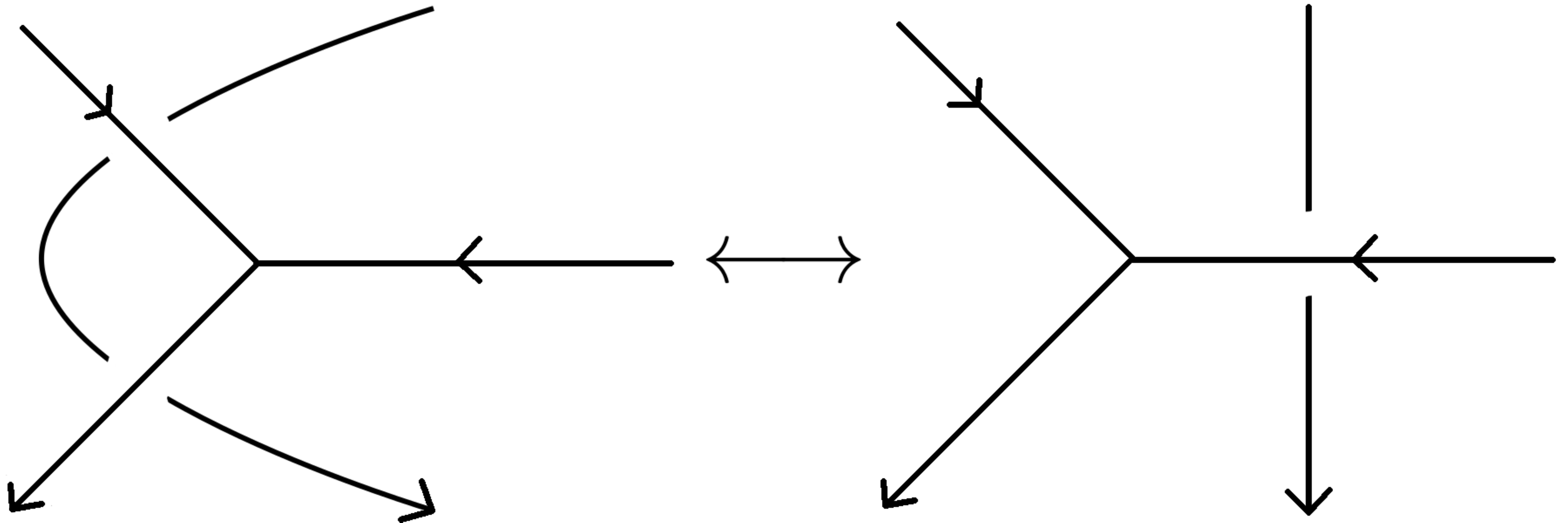} & 
\includegraphics [height=0.4in]{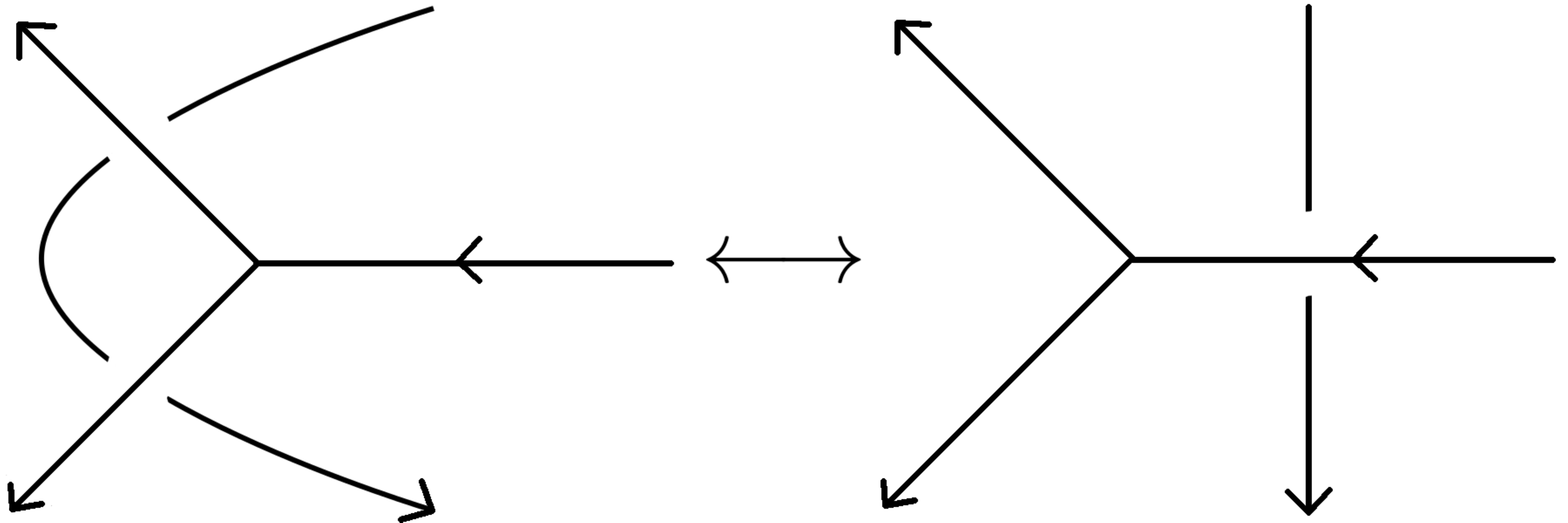}\\
\vspace{1em} \\
$R5.17$ & $R5.18$ & $R5.19$ & $R5.20$ \\
\includegraphics[height=0.4in]{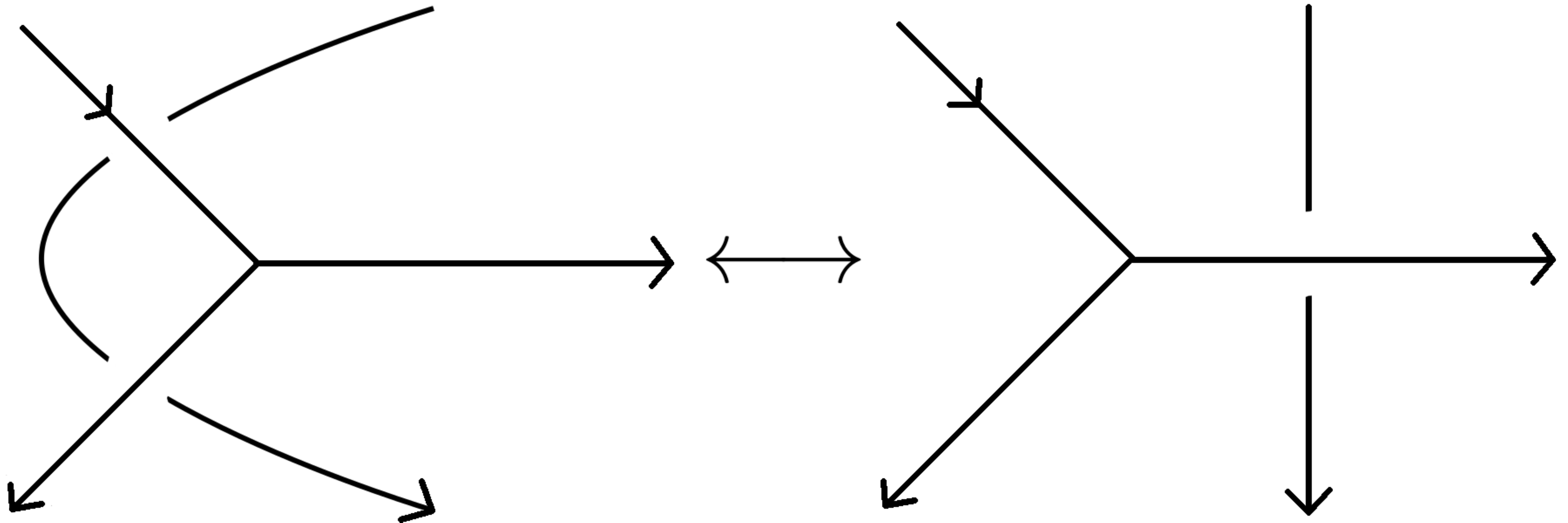} & 
\includegraphics[height=0.4in]{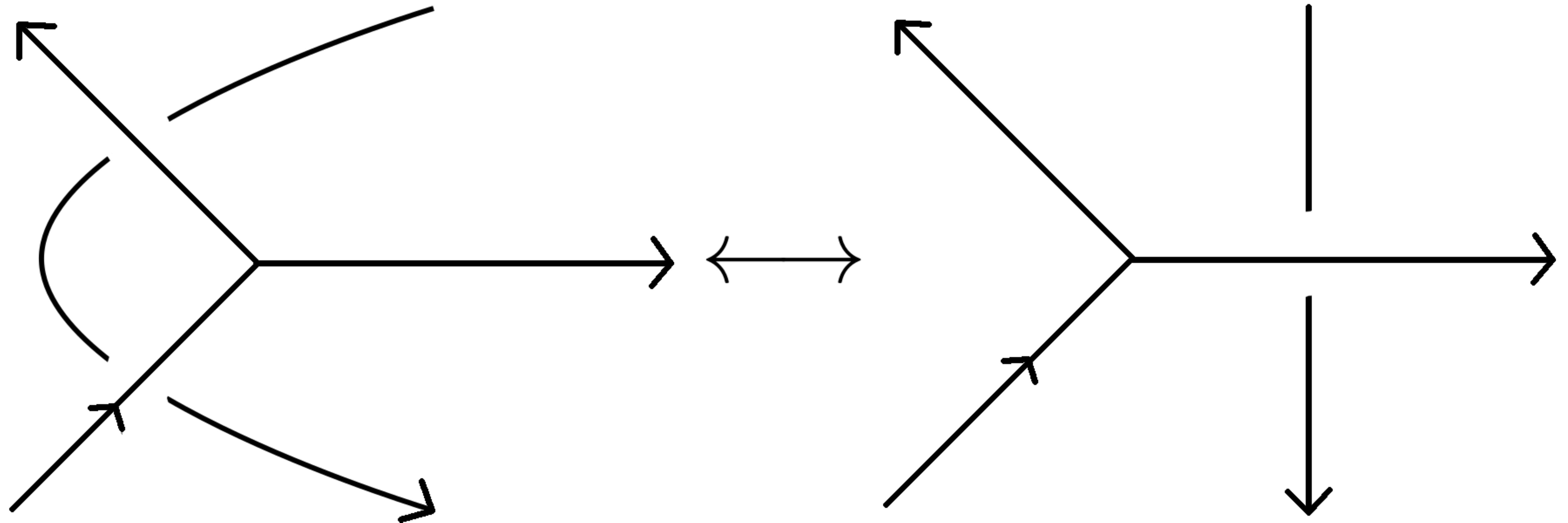} &
\includegraphics[height=0.4in]{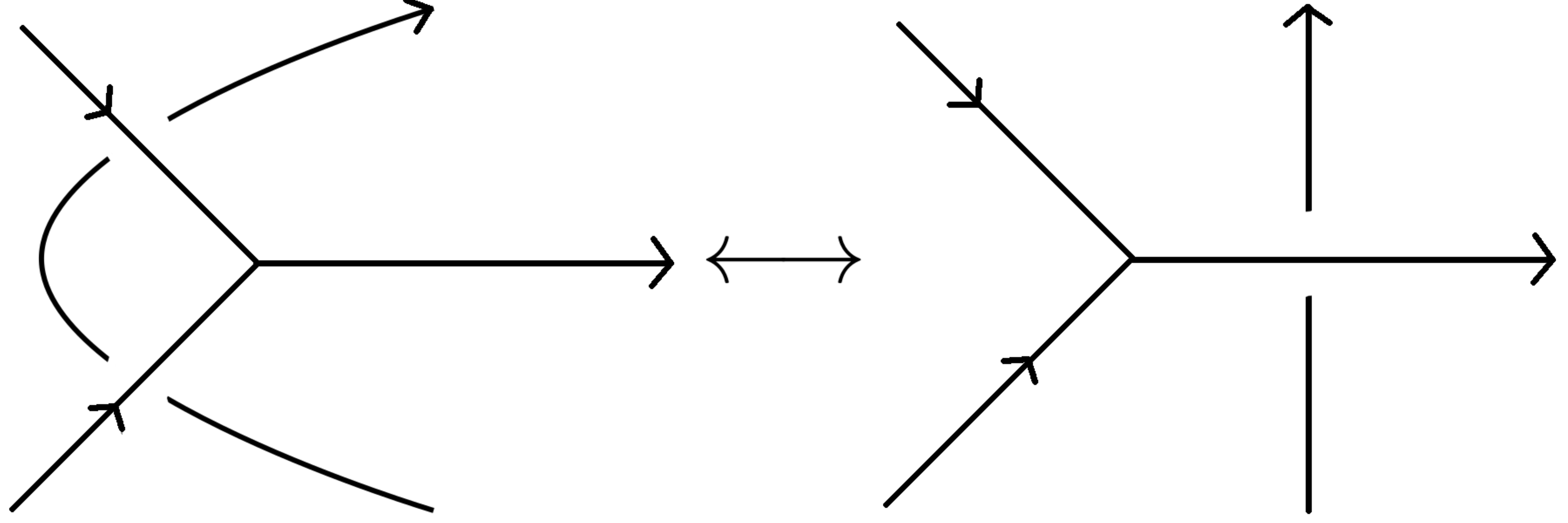}&
\includegraphics[height=0.4in]{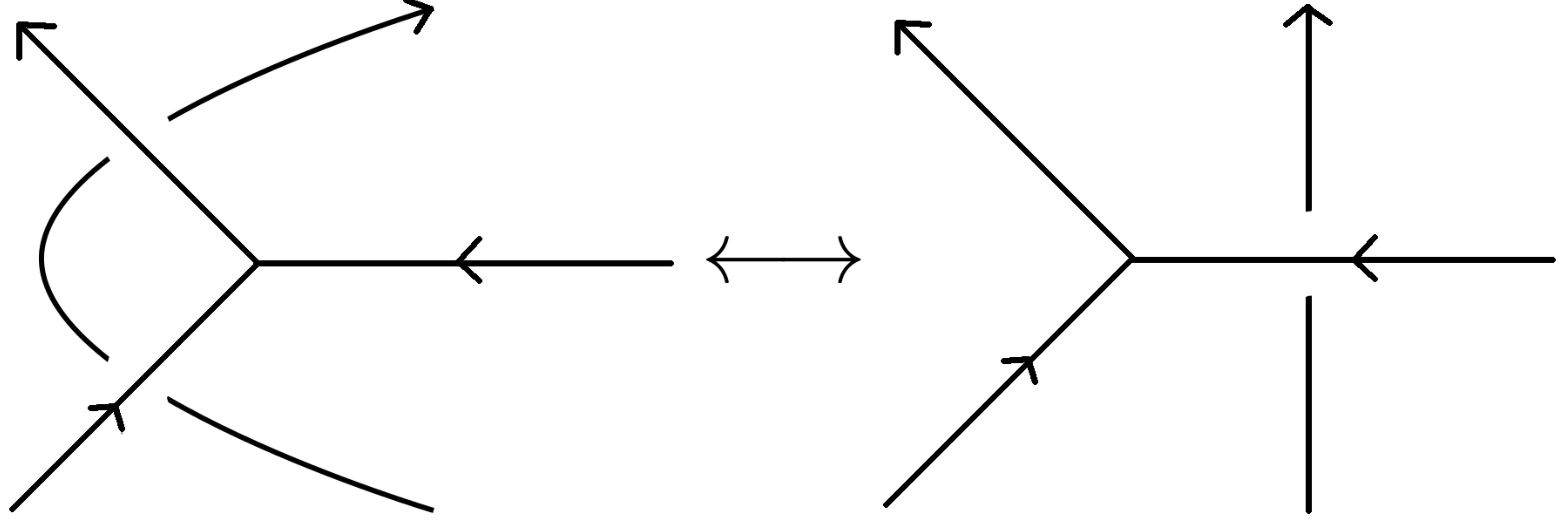} \\
\vspace{1em} \\
$R5.21$ & $R5.22$ & $R5.23$ & $R5.24$ \\
\includegraphics[height=0.4in]{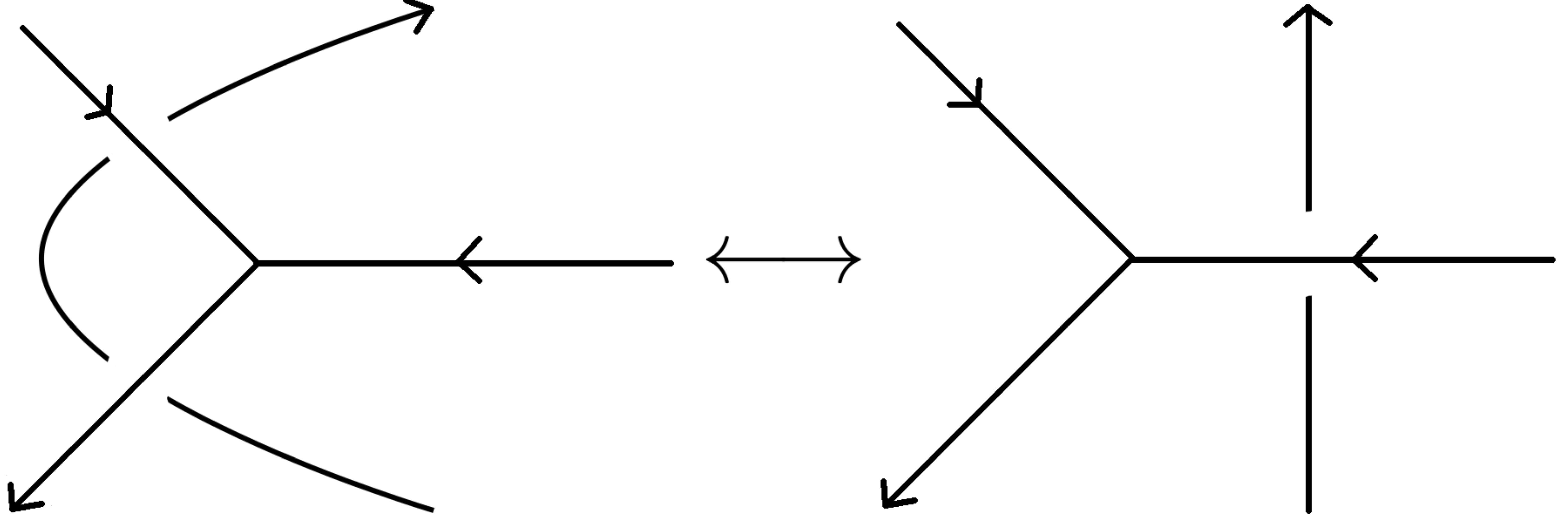}&
\includegraphics[height=0.4in]{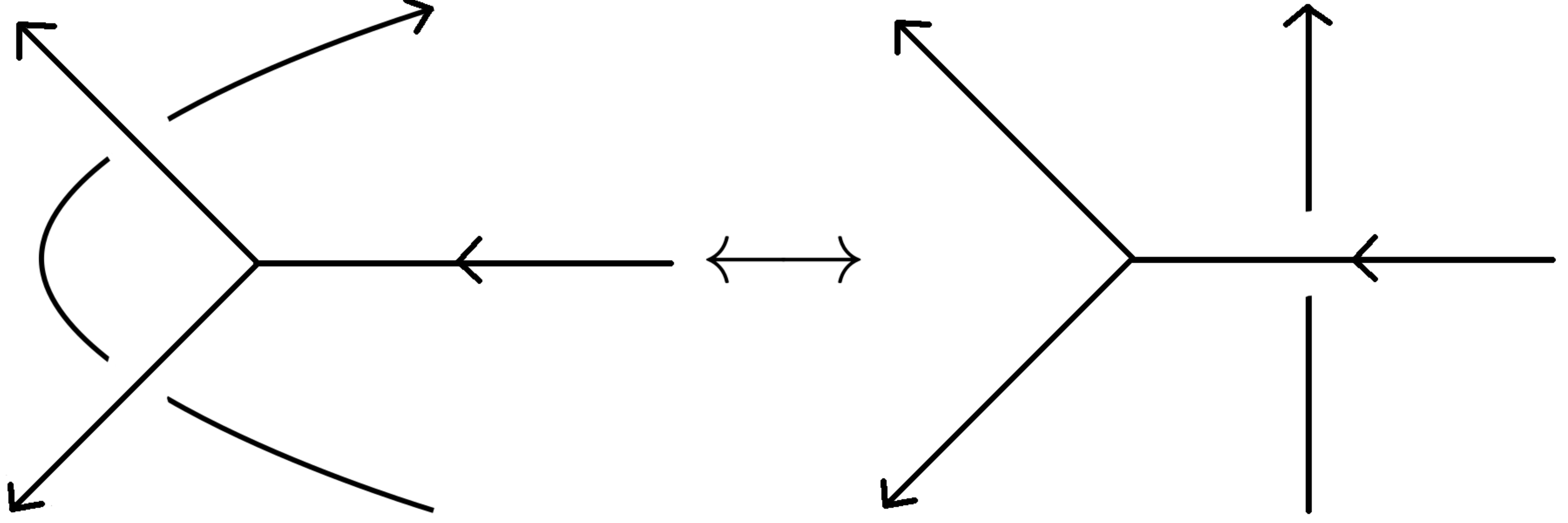} &
\includegraphics[height=0.4in]{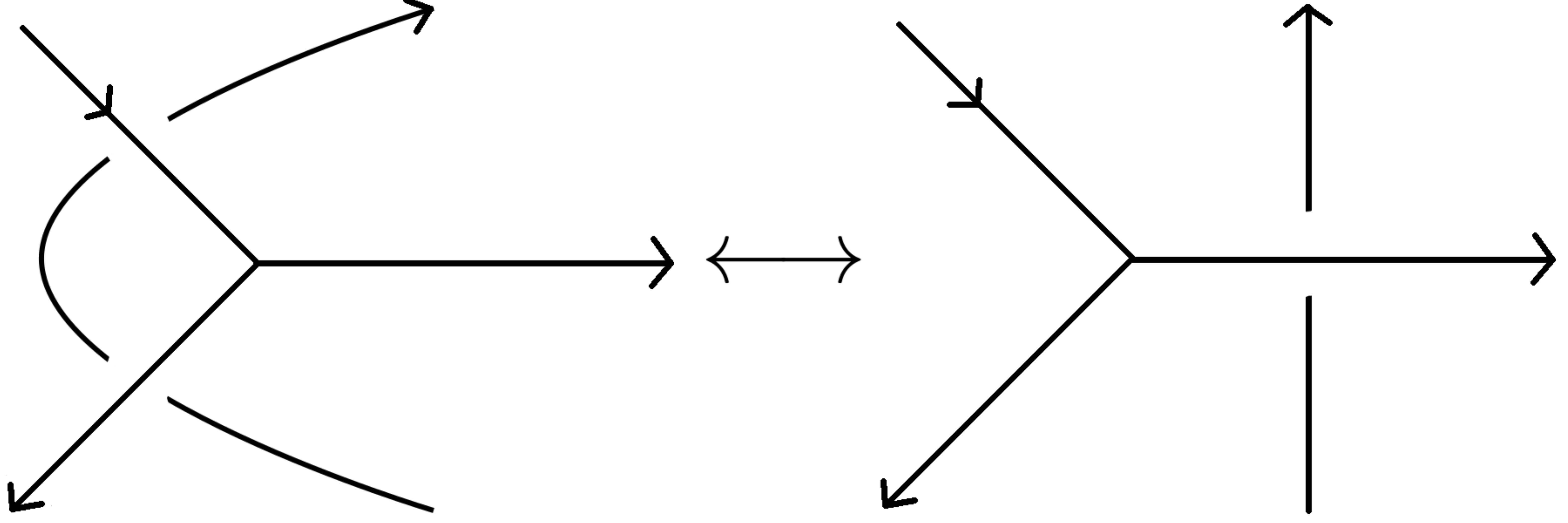} &
\includegraphics[height=0.4in]{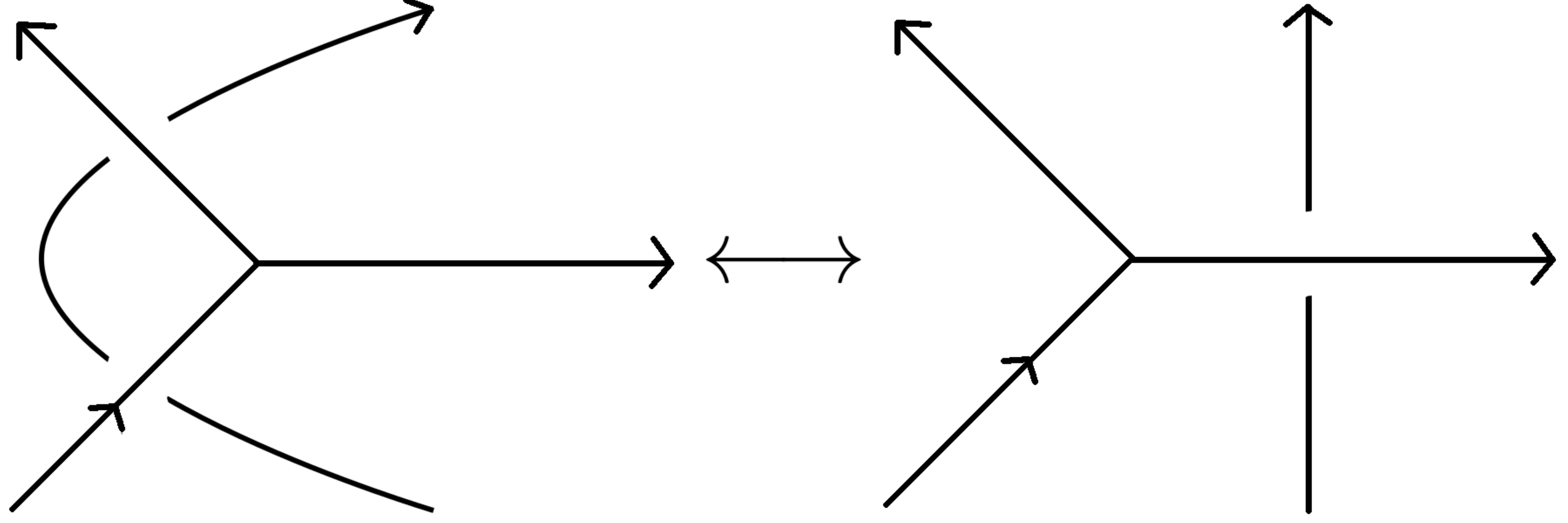}
\end{tabular}
\end{center}

\begin{lemma} We can obtain a generating set of Reidemeister moves on trivalent spacial graphs by taking one move from each of the following sets
\begin{center}
\begin{tabular}{l}
\{R4.1, R4.2, R4.3, R4.7, R4.8, R4.9\},\\
\{R4.4, R4.5, R4.6, R4.10, R4.11, R4.12\} \\
\\
\{R5.1, R5.2, R5.3, R5.7, R5.8, R5.9\} \\
\{R5.4, R5.5, R5.6, R5.10, R5.11, R5.12\} \\
\{R5.13, R5.14, R5.15, R5.19, R5.20, R5.21\} \\
\{R5.16, R5.17, R5.18, R5.22, R5.23, R5.24\} \\
\end{tabular}
\end{center}
together with any minimal generating set of oriented classical Reidemeister moves.
\end{lemma}

For a list of  minimal generating sets of oriented classical Reidemeister moves, see \cite{P}.

\begin{proof}
Through the sequence of moves given by

\begin{center}
\includegraphics[height=0.8in]{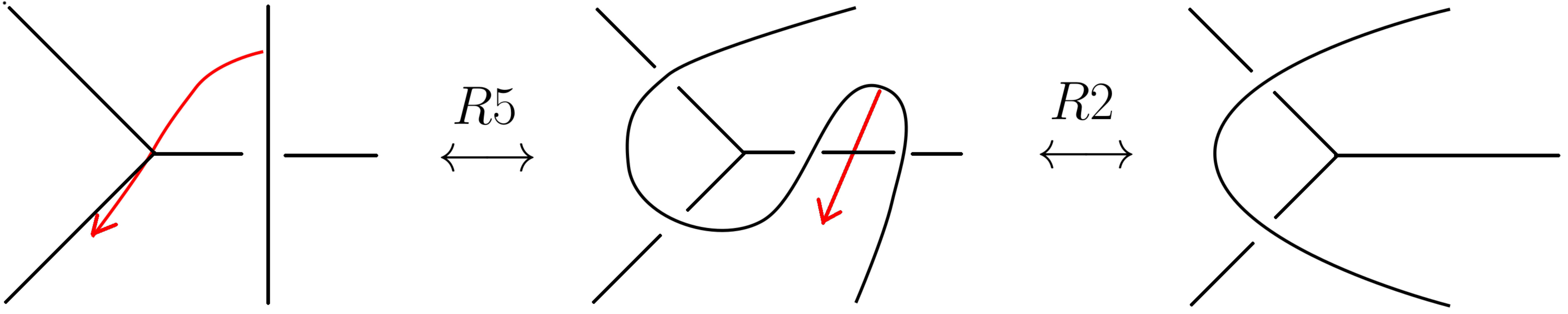} 
\end{center}
and by assigning different admissible combinations of orientations to each of the edges, we see that, together with $R2$,  we have the sequences of implications
$$ R5.1 \Rightarrow R5.9 \Rightarrow R5.2 \Rightarrow R5.7 \Rightarrow R5.3 \Rightarrow R5.8 \Rightarrow R5.1$$
and
$$R5.4 \Rightarrow R5.12 \Rightarrow R5.5 \Rightarrow R5.10 \Rightarrow R5.6 \Rightarrow R5.11 \Rightarrow R5.4$$
so the sets
$$\{R5.1, R5.2, R5.3, R5.7, R5.8, R5.9\}$$
and
$$\{R5.4, R5.5, R5.6, R5.10, R5.11, R5.12\}$$
give sets of equivalent moves.

Similarly, using the sequence of moves

\begin{center}
\includegraphics[height=0.8in]{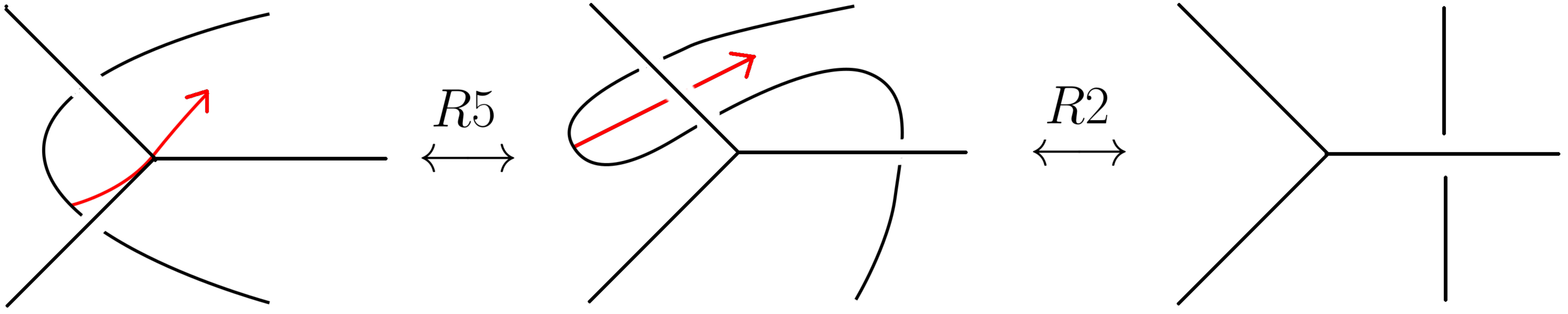}
\end{center}
and assigning different admissible orientations to each of the edges, we see that, together with $R2$, we have the sequences of implications
$$R5.13 \Rightarrow R5.21 \Rightarrow R5.14 \Rightarrow R5.19 \Rightarrow R5.15 \Rightarrow R5.20\Rightarrow R5.13$$
and
$$R5.16 \Rightarrow R5.24 \Rightarrow R5.17 \Rightarrow R5.22 \Rightarrow R5.18 \Rightarrow R5.23 \Rightarrow R5.16$$
Thus the sets
$$\{R5.13, R5.14, R5.15, R5.19, R5.20, R5.21\}$$
and
$$\{R5.16, R5.17, R5.18, R5.22, R5.23, R5.24\}$$
give sets of equivalent moves.

Now, we see that

\begin{center}
\includegraphics[height=0.8in]{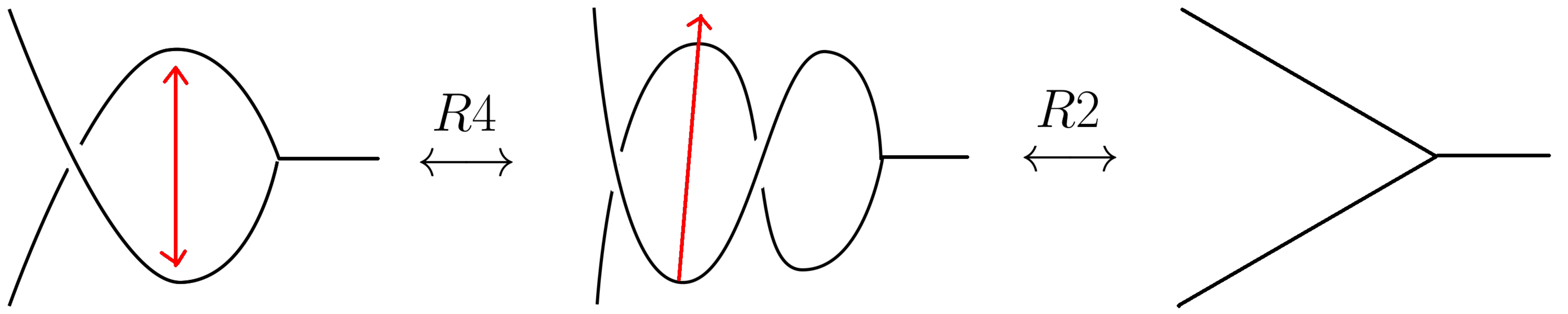}
\end{center}
tells us that
\begin{center}
\begin{tabular}{l@{\hspace{1in}}l@{\hspace{1in}}l}
$R4.1 \Rightarrow R4.7$ & $R4.2 \Rightarrow R4.9$ & $R4.3 \Rightarrow R4.8$ \\
$R4.4 \Rightarrow R4.10$ & $R4.5 \Rightarrow R4.12$ & $R4.6 \Rightarrow R4.11$
\end{tabular}
\end{center}
and

\begin{center}
\includegraphics[height=0.8in]{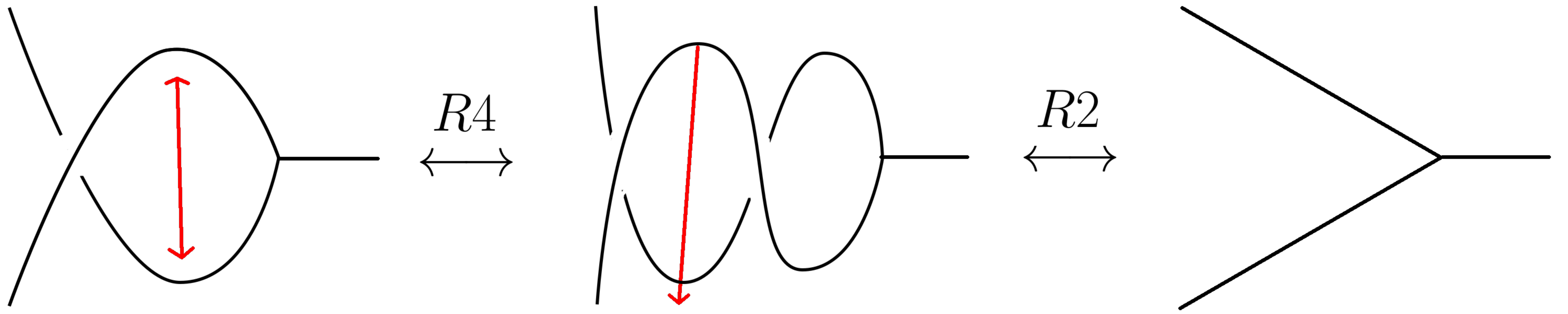}
\end{center}
tells us 
\begin{center}
\begin{tabular}{l@{\hspace{1in}}l@{\hspace{1in}}l}
$R4.7 \Rightarrow R4.1$ & $R4.8 \Rightarrow R4.3$ & $R4.9 \Rightarrow R4.2$ \\
$R4.10 \Rightarrow R4.4$ & $R4.11 \Rightarrow R4.6$ & $R4.12 \Rightarrow R4.5$
\end{tabular}
\end{center}
 so the sets
 $$\{R4.1, R4.7\}, \, \{R4.2, R4.9\}, \, \{R4.3, R4.8\}, \,\{R4.4, R4.10\},\, \{R4.5, R4.12\},\, \text{and }\{R4.6, R4.11\}$$
 are sets of equivalent moves, if we assume $R2$.
 
 Furthermore, if we assume the $R5$-type moves and $R1$, the sequence of moves
 
 \begin{center}
 \includegraphics[height=0.8in]{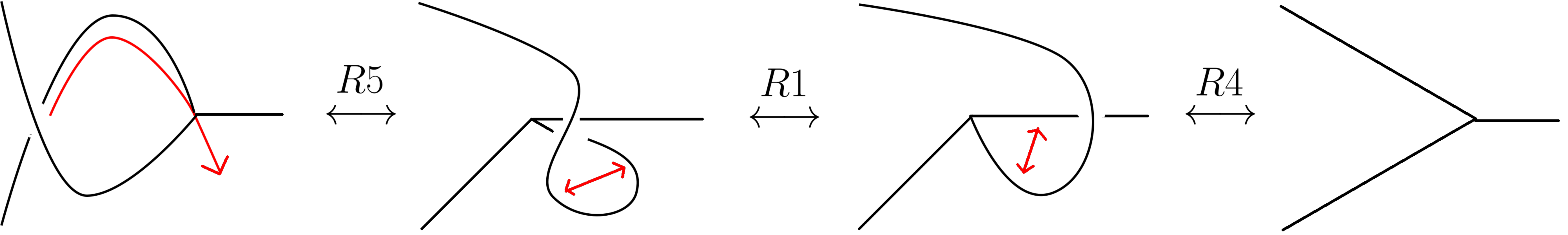}
 \end{center}
gives us that 
\begin{center}
\begin{tabular}{l@{\hspace{0.7in}}l@{\hspace{0.7in}}l@{\hspace{0.7in}}l}
$R4.2 \Rightarrow R4.7$ & $R4.1 \Rightarrow R4.9$ & $R4.6 \Rightarrow R4.10$ & $R4.4 \Rightarrow R4.11$
\end{tabular}
\end{center}
 
 Similarly, we have
 
 \begin{center}
 \includegraphics[height=0.8in]{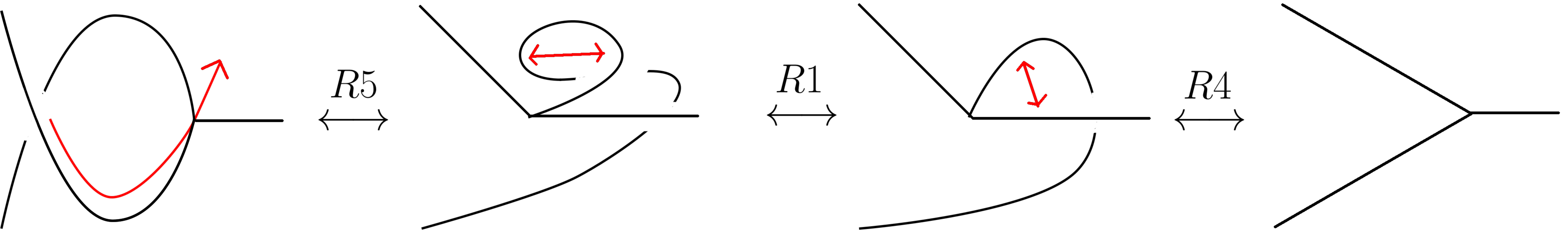}
 \end{center}
 
 which tells us that, together with the classical Reidemeister moves and the $R5$-types moves, 
 \begin{center}
\begin{tabular}{l@{\hspace{0.7in}}l@{\hspace{0.7in}}l@{\hspace{0.7in}}l}
$R4.3 \Rightarrow R4.7$ & $R4.1 \Rightarrow R4.8$ & $R4.5 \Rightarrow R4.10$ & $R4.4 \Rightarrow R4.12$
\end{tabular}
\end{center}
So, we have sets of equivalent moves 
 $$\{R4.1, R4.2, R4.3, R4.7, R4.8, R4.9\}$$
 and 
 $$\{R4.4, R4.5, R4.6, R4.10, R4.11, R4.12\}.$$

Putting this all together, we have equivalence classes of moves given by

\begin{center}
\begin{tabular}{l}
\{R4.1, R4.2, R4.3, R4.7, R4.8, R4.9\},\\
\{R4.4, R4.5, R4.6, R4.10, R4.11, R4.12\} \\
\\
\{R5.1, R5.2, R5.3, R5.7, R5.8, R5.9\} \\
\{R5.4, R5.5, R5.6, R5.10, R5.11, R5.12\} \\
\{R5.13, R5.14, R5.15, R5.19, R5.20, R5.21\} \\
\{R5.16, R5.17, R5.18, R5.22, R5.23, R5.24\} \\
\end{tabular}
\end{center}

So, for a reduced Reidemeister moves set, we only need to choose one from each class, together with a minimal generating set of classical oriented Reidemeister moves.
\end{proof}

By choosing a specific move from each equivalence class, we have the following corollary.

\begin{corollary}
The  moves

\begin{center}
\begin{tabular}{ccc}
$R4.1$ & \;\;\; &$R4.10$\\
\includegraphics[height=0.8in]{rreduce-r4-1} &&
\includegraphics[height=0.8in]{rreduce-r4-10}\\
\vspace{1em} \\
$R5.7$ && $R5.10$ \\
\includegraphics[height=0.8in]{rreduce-r5-7} &&
\includegraphics[height=0.8in]{rreduce-r5-10}\\
\vspace{1em} \\
$R5.13$ && $R5.16$ \\
\includegraphics[height=0.8in]{rreduce-r5-13} &&
\includegraphics[height=0.8in]{rreduce-r5-16}\\
\end{tabular}
\end{center}
together with a generating set of oriented classical Reidemeister moves
form a generating set of oriented Reidemeister moves for trivalent spatial graphs.

\end{corollary}

\section{\large \textbf{Niebrzydowski Algebras}}\label{B}

We begin with a modification of a definition from \cite{MN}.

\begin{definition}
Let $X$ be a set. A ternary operation $[,,]:X\times X\times X\to X$ on $X$
is a \textit{Niebrzydowski tribracket} or just a \textit{tribracket}
if it satisfies the conditions
\begin{itemize}
\item[(i)] Any three of the four $\{a,b,c,d\}$ in the equation $[a,b,c]=d$
determines the fourth,
\item[(ii)] For all $a,b,c,d\in X$ we have
\[\begin{array}{rcl}
{[a,b,[b,c,d]]} & = & [a,[a,b,c],[[a,b,c],c,d]]\\
{[[a,b,c],c,d]} & = & [[a,b,[b,c,d]],[b,c,d],d].
\end{array}\]
\end{itemize} 
\end{definition}

\begin{example}
Any module over a ring $R$ has a Niebrzydowski bracket structure defined
by setting 
\[[a,b,c,]=xa-xyb+yc\]
where $x,y\in R^{\times}$ are units in $R$. This structure is known as an
\textit{Alexander tribracket}; see \cite{NP} for more.
\end{example}

\begin{example}\label{ex:1}
We can specify a Niebrzydowski algebra structure on a finite set 
$X=\{1,2,\dots, n\}$ by listing the operation tables. Specifically, the
Niebrzydowski bracket can be expressed as a $3$-tensor, i.e. a vector
whose entries are matrices, with the convention that to find $[a,b,c]$
we look in matrix $a$, row $b$ column $c$. For example, in the tribracket
structure on $\{1,2,3\}$ specified by
\[
\left[\left[\begin{array}{rrr}
1& 2& 3 \\
3& 1& 2 \\
2& 3& 1
\end{array}\right],
\left[\begin{array}{rrr}
2& 3& 1 \\
1& 2& 3 \\
3& 1& 2
\end{array}\right],
\left[\begin{array}{rrr}
3& 1& 2 \\
2& 3& 1 \\
1& 2& 3
\end{array}\right]\right]
\]
we have $[1,2,3]=2$ and $[2,3,1]=3$.
\end{example}

\begin{definition}
Let $X$ be a set with a Niebrzydowski tribracket $[,,]$. Then $X$ is a
\textit{Niebrzydowski algebra} if $X$ has a partially defined product 
$a,b\mapsto ab$ satisfying
\begin{itemize}
\item[(i)] Any two of the three $\{a,b,c\}$ in $ab=c$ determines the third,
\item[(ii)] For all $a,b\in X$ we have
\[[a,ab,b]=ab \quad (iv.i)\]
and
\item[(iii)] For all $a,b,c\in X$ we have
\[\begin{array}{rcll}
a[a,b,c] & = & [a,b,bc] & (v.i)\\
{[a,b,c]c} & = & [ab,b,c] & (v.ii)\\
{[a,b,c]} & = & [[a,b,bc],bc,c] & (v.iii)\\
{[a,b,c]} & = & [a,ab,[ab,b,c]] & (v.iv). 
\end{array}
\]
\end{itemize}
\end{definition}

\begin{example}
The tribracket in example \ref{ex:1} has eight partial products, given 
by the operation matrices
\[
\left[\begin{array}{rrr}1& 3 & 2 \\ 3 & 2 & 1 \\ 2 & 1 & 3 \end{array}\right],
\left[\begin{array}{rrr}1& 3 & - \\ - & 2 & 1 \\ 2 & - & 3 \end{array}\right],
\left[\begin{array}{rrr}1& - & 2 \\ 3 & 2 & - \\ - & 1 & 3 \end{array}\right],
\left[\begin{array}{rrr}-& 3 & 2 \\ 3 & - & 1 \\ 2 & 1 & - \end{array}\right],
\]\[
\left[\begin{array}{rrr}1& - & - \\ - & 2 & - \\ - & - & 3 \end{array}\right],
\left[\begin{array}{rrr}-& 3 & - \\ - & - & 1 \\ 2 & - & - \end{array}\right],
\left[\begin{array}{rrr}-& - & 2 \\ 3 & - & - \\ - & 1 & - \end{array}\right],
\left[\begin{array}{rrr}-& - & - \\ - & - & - \\ - & - & - \end{array}\right]
\]
where a $-$ indicates an undefined product.
\end{example}

These elements of $X$ are interpreted as colors for regions in the planar 
complement of a $Y$-oriented trivalent spatial graph diagram according to 
the rule
\[\includegraphics{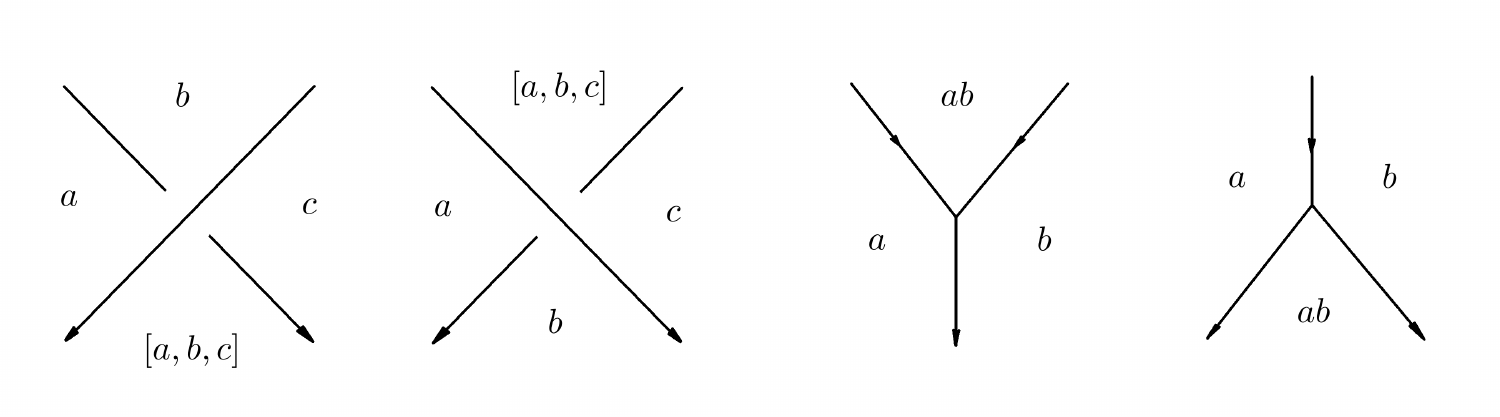}\]

The axioms come from our preferred generating set of $Y$-oriented
Reidemeister moves. From the R4 moves, we have
\[\includegraphics{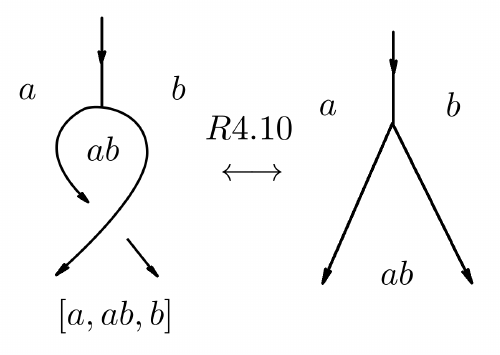}\quad \includegraphics{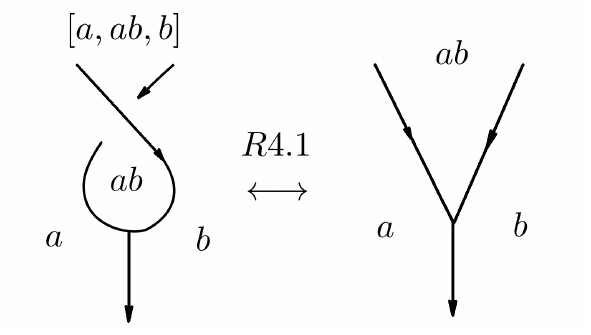}\]
and from the R5 moves we have
\[\includegraphics{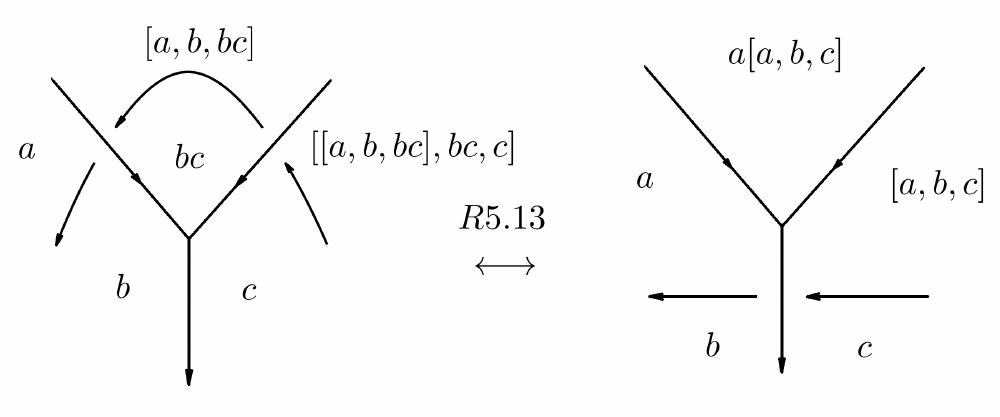}\]
\[\includegraphics{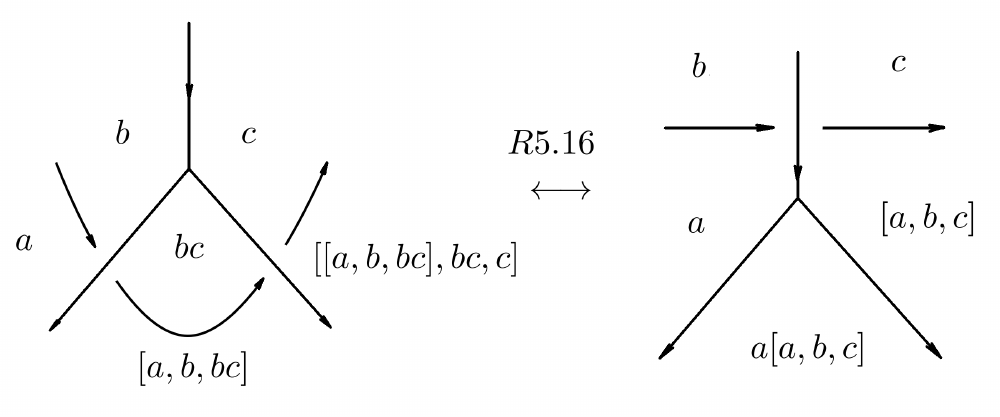}\]
\[\includegraphics{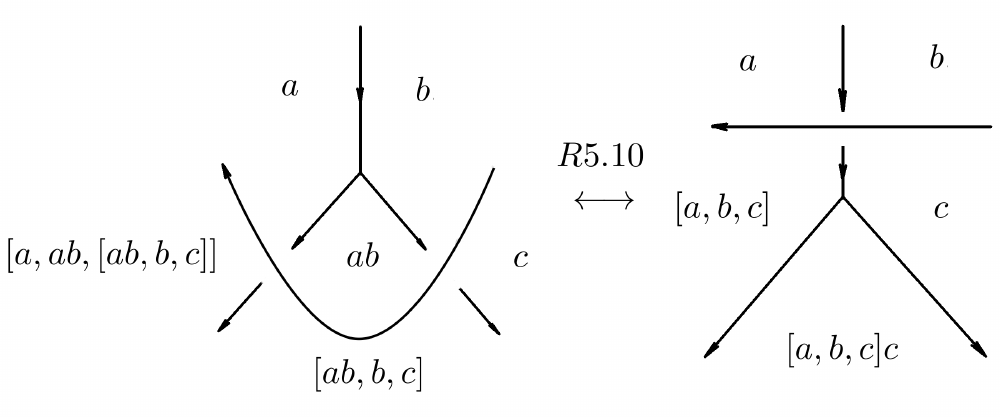}\]
and
\[\includegraphics{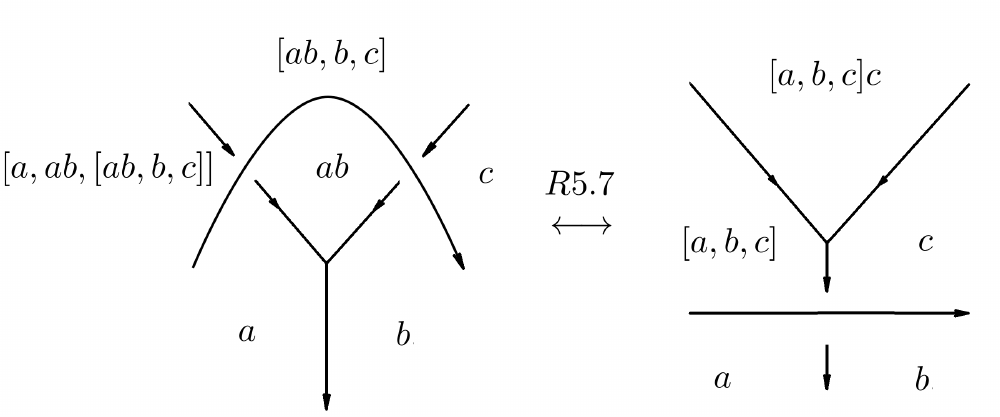}.\]

\begin{remark}
The various oriented IH moves impose the idempotency-style condition that 
\[a(ab)=ab=(ab)b\]
for all $a,b\in X$. We illustrate with two of the possible $Y$-oriented
IH-moves; the other cases are similar.
\[\includegraphics{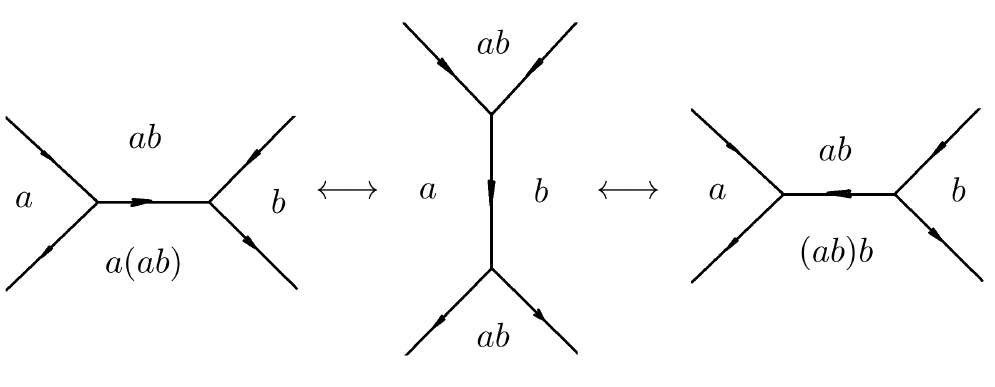}\]
The left and right invertibility of the partial product then require that 
$ab=a=b$, so
this condition implies that the multiplication is only defined 
for equal operands,
i.e., if $a\ne b$ then $ab$ is not defined, and that we must have $aa=a$.
Thus, for  handlebody-links the rule is that all three region 
colors around a vertex must be equal. We will call such a Niebrzydowski algebra
an \textit{idempotent Niebrzydowski algebra}.
\end{remark}

Then by construction, we have
\begin{theorem}
Let $X$ be a Niebrzydowski algebra and $\Gamma$ a $Y$-oriented trivalent
spatial graph. Then the number $\Phi_X^{\mathbb{Z}}(\Gamma)$ of $X$-colorings
of any diagram of $\Gamma$ is invariant under $Y$-oriented Reidemeister
moves. If $X$ is idempotent, then $\Phi_X^{\mathbb{Z}}(\Gamma)$ is an
invariant of $Y$-oriented handlebody-links.
\end{theorem}

\begin{example}\label{ex:2}
The $Y$-oriented unknotted theta graph below has nine colorings by 
the Niebrzydowski algebra as depicted.
\[X=\left[\left[\begin{array}{rrr}
1& 2& 3 \\
3& 1& 2 \\
2& 3& 1
\end{array}\right],
\left[\begin{array}{rrr}
2& 3& 1 \\
1& 2& 3 \\
3& 1& 2
\end{array}\right],
\left[\begin{array}{rrr}
3& 1& 2 \\
2& 3& 1 \\
1& 2& 3
\end{array}\right]\right], 
\left[\begin{array}{rrr}
1& 3& 2 \\
3& 2& 1 \\
2& 1& 3
\end{array}\right]
\]
\[\scalebox{0.95}{\includegraphics{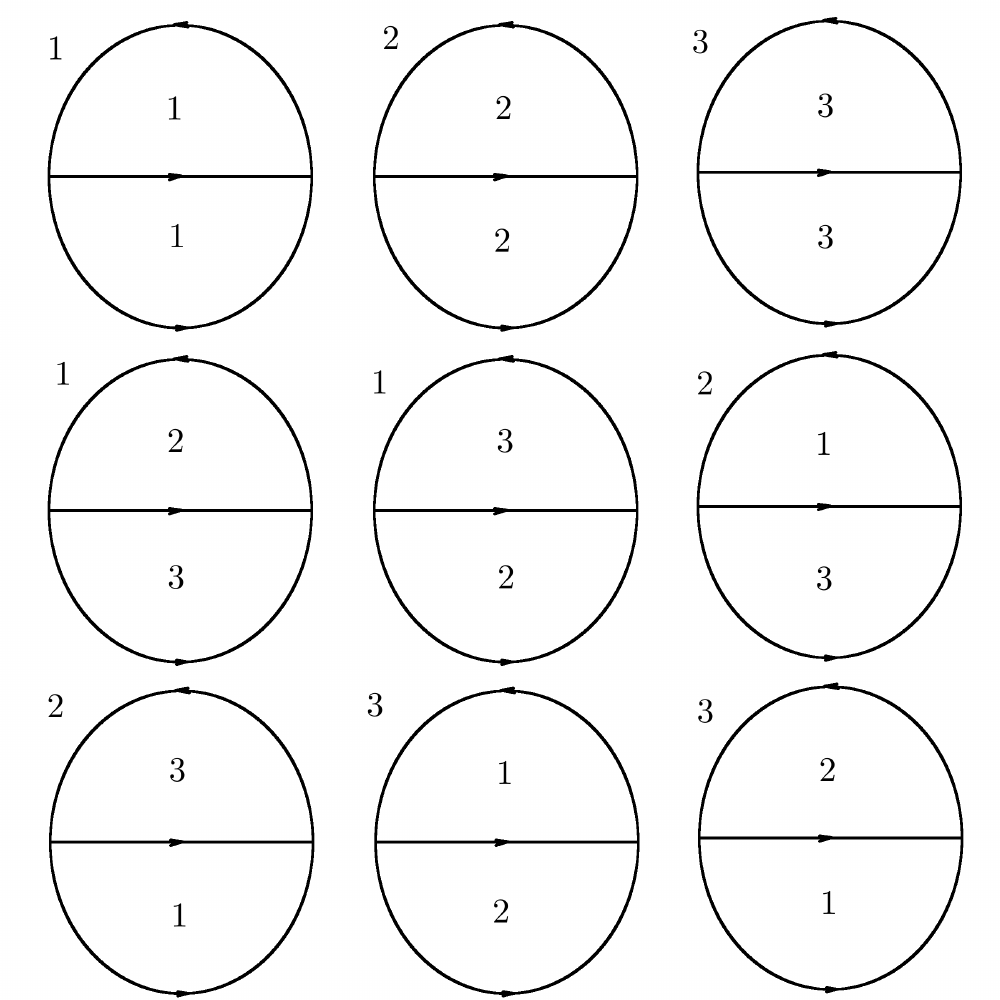}}\]
Hence we have counting invariant value $\Phi_X^{\mathbb{Z}}(\Gamma)=9$.
This is distinct from the $Y$-oriented handcuff graph, which has only
three $X$-colorings. 
\[\includegraphics{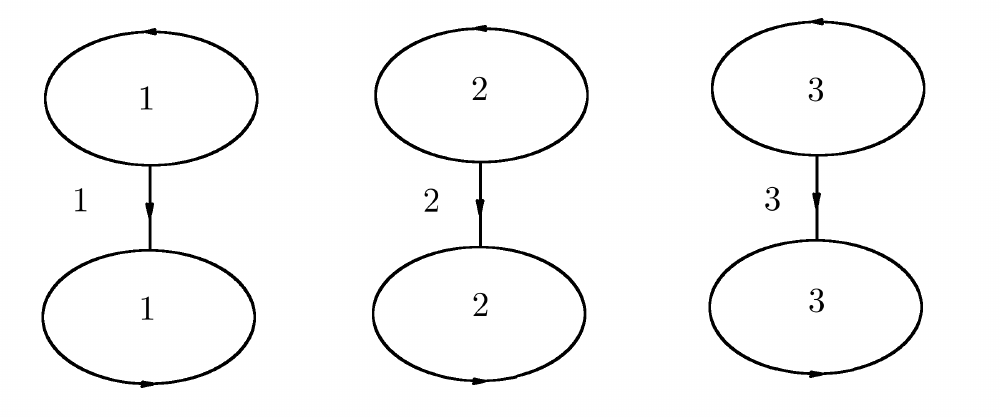}\]
We note that while these two graphs represent
the same handlebody-knot, $X$ is not an idempotent Niebrzydowski
algebra and thus can distinguish different trivalent spatial graphs which 
represent the same handlebody-knot. If we instead use the partial 
product given by the matrix
\[\left[\begin{array}{rrr}
1& -& - \\
-& 2& - \\
-& -& 3
\end{array}\right]\]
then the resulting Niebrzydowski algebra is idempotent and both graphs
have only the three constant colorings, reflecting the fact that they represent
the same handlebody-knot.
\end{example}

\begin{example}\label{ex:hbdy}
Idempotent Niebrzydowski algebra counting invariants can detect when handlebody-links have different genera. For example, the diagram on the left
\[\includegraphics{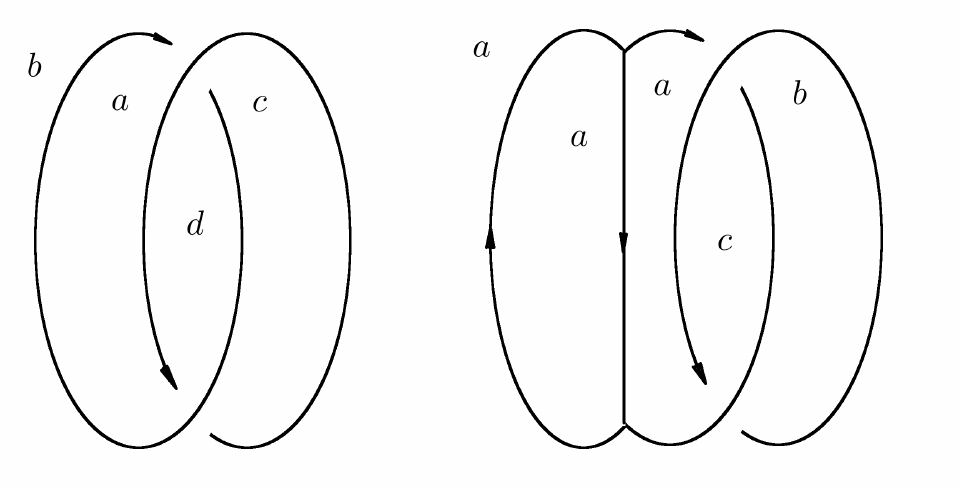}\]
represents a handlebody-link of two genus 1 components (i.e., a classical link) 
while the one on the right represents a handlebody-link of a genus 2 component 
with a genus 1 component. These are distinguished by the counting invariant 
with respect to the idempotent Niebrzydowski algebra 
\[X=\left[\left[\begin{array}{rrr}
1& 2& 3 \\
3& 1& 2 \\
2& 3& 1
\end{array}\right],
\left[\begin{array}{rrr}
2& 3& 1 \\
1& 2& 3 \\
3& 1& 2
\end{array}\right],
\left[\begin{array}{rrr}
3& 1& 2 \\
2& 3& 1 \\
1& 2& 3
\end{array}\right]\right], 
\left[\begin{array}{rrr}
1& 0& 0 \\
0& 2& 0 \\
0& 0& 3
\end{array}\right]\]
with counting invariant values $27$ and $3$ respectively: checking the 
operation tensor of $X$ we verify that $[a,b,c]=d$ iff $[a,d,c]=b$, so the 
regions labeled $a,b,c$ can be freely chosen and we have $3^3=27$ $X$-colorings
of the Hopf link, while for the handlebody-link on the right the coloring 
conditions $[a,a,b]=c$ and $[a,c,b]=a$ are satisfied only for $a=b=c$, yielding
only three colorings.
\end{example}

\begin{example}
Let $X$ be the Niebrzydowski algebra given by
\[X=\left[\left[\begin{array}{rrr}
1& 2& 3 \\
3& 1& 2 \\
2& 3& 1
\end{array}\right],
\left[\begin{array}{rrr}
2& 3& 1 \\
1& 2& 3 \\
3& 1& 2
\end{array}\right],
\left[\begin{array}{rrr}
3& 1& 2 \\
2& 3& 1 \\
1& 2& 3
\end{array}\right]\right], 
\left[\begin{array}{rrr}
-& 3& - \\
-& -& 1 \\
2& -& -
\end{array}\right]
\]
and consider the two trivalent spatial graphs below.
\[\includegraphics{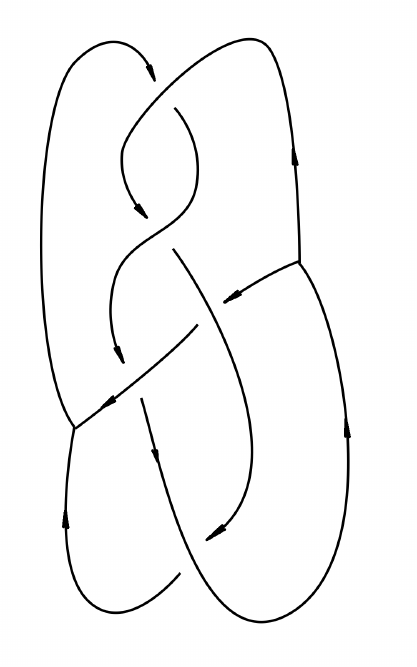}
\quad \includegraphics{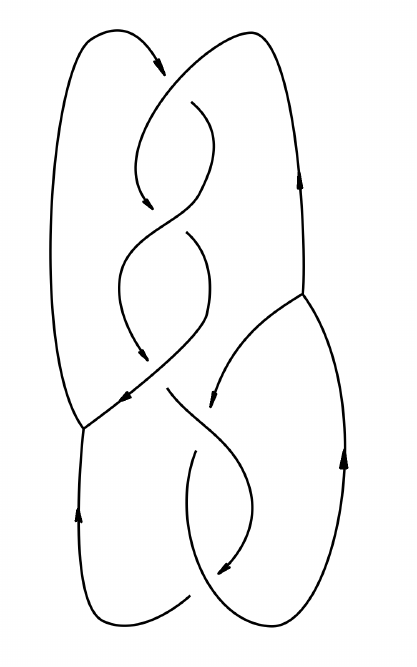}
\]
The spatial graph $K_1$ on the left has three $X$-colorings
\[\includegraphics{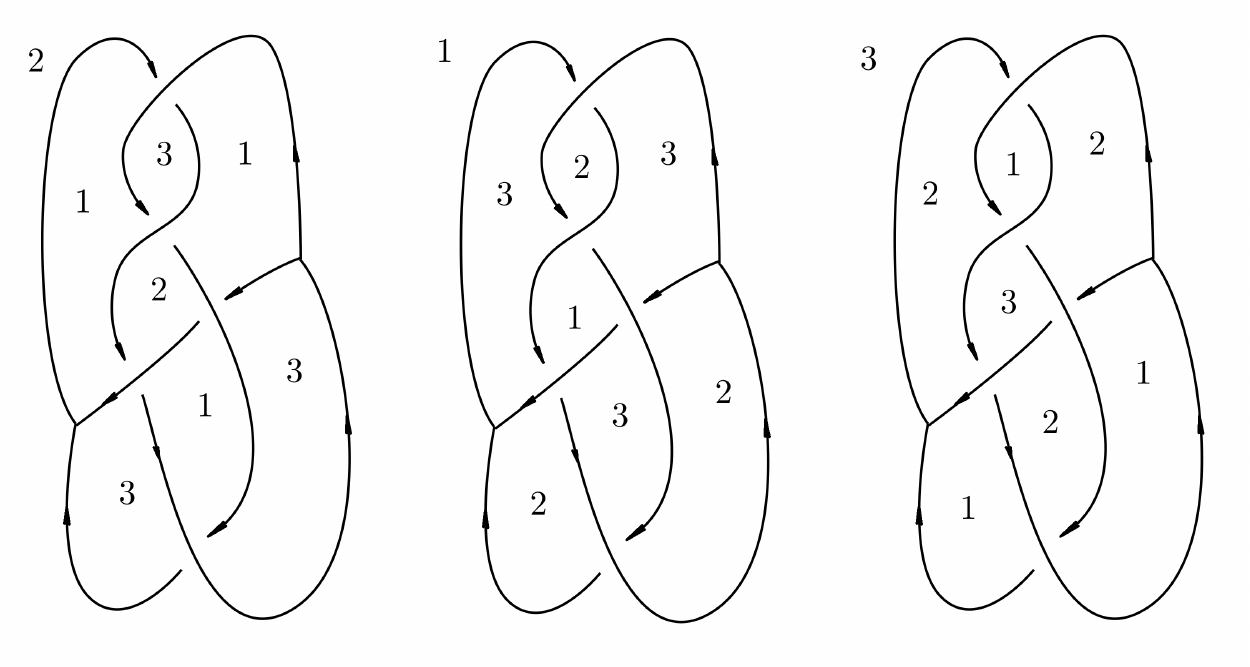}\]
for a counting invariant value of $\Phi_X^{\mathbb{Z}}(K_1)=3$. 
For the spatial graph $K_2$ on the right, we can observe that the partial
multiplication structure of $X$ means that the only possible colorings
around the vertices are as depicted with 
$(a,b,c)\in\{(1,2,3),(2,3,1),(3,1,2)\}$.
\[\includegraphics{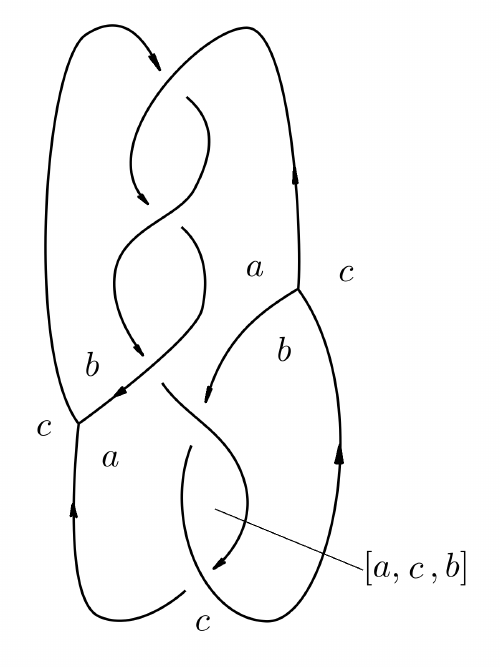}\]
However, for such a coloring to be valid, we would need
\[[a,[a,c,b],b]=a\]
and we verify that none of the three possibilities works:
\[\begin{array}{rcl}
{[}1,[1,3,2],2] & = & [1,3,2] = 3\ne 1 \\
{[}2,[2,1,3],3] & = & [2,1,3]=1\ne 2 \\
{[}3,[3,2,1],1] & = & [3,2,1]=2\ne 3.
\end{array}\]
Hence there are no $X$-colorings of $K_2$ and the counting invariant
distinguishes these spatial graphs.
\end{example}

\begin{example}
Using our \texttt{Python} code, we find that the trivalent spatial graphs 
\[\includegraphics{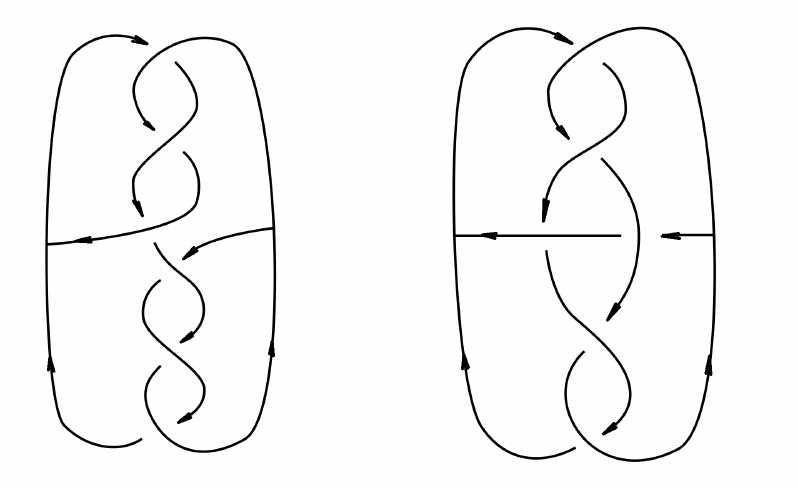}\]
are distinguished by their numbers of colorings by the Niebrzydowski algebra
$X$ specified by
\[
\left[
\left[\begin{array}{rrrr}
1& 2& 3& 4\\
4& 1& 2& 3\\
3& 4& 1& 2\\
2& 3& 4& 1
\end{array}\right],\left[\begin{array}{rrrr}
2& 3& 4& 1\\
1& 2& 3& 4\\
4& 1& 2& 3\\
3& 4& 1& 2
\end{array}\right],\left[\begin{array}{rrrr}
3& 4& 1& 2\\
2& 3& 4& 1\\
1& 2& 3& 4\\
4& 1& 2& 3
\end{array}\right],\left[\begin{array}{rrrr}
4& 1& 2& 3\\
3& 4& 1& 2\\
2& 3& 4& 1\\
1& 2& 3& 4
\end{array}\right]\right],
\left[\begin{array}{rrrr}
1& -& 2& -\\
-& 2& -& 3\\
4& -& 3& -\\
-& 1& -& 4
\end{array}\right]
\]
with eight and four colorings respectively for the spatial graphs on the 
left and on the right.
\end{example}

\section{\large \textbf{Questions}}\label{Q}

We end with a few questions for future research.

We can't help but notice that the fully-defined multiplication table for the
tribracket in example \ref{ex:2} is a quandle table, namely the dihedral 
quandle of three elements. Is this a coincidence?

Because of the idempotency requirement, handlebody-link colorings tend to be
fairly monochromatic. How many crossings in a handlebody knot (single component)
are needed before we get a non-monochromatic coloring by an idempotent 
Niebrzydowski algebra?

As always, what enhancements of this counting invariant are possible? What 
are the connections with qualgebra colorings and qualgebra homology?

\bibliography{pg-sn-st}{}
\bibliographystyle{abbrv}

\bigskip

\noindent\textsc{Department of Mathematics \\ 
Claremont McKenna College \\
850 Columbia Ave. \\
Claremont, CA 91711}

\bigskip
\noindent\textsc{Department of Mathematics \\ 
South Hall, Room 6607\\
University of California \\
Santa Barbara, CA 93106}

\end{document}